%% file: main_arxiv.tex
\documentclass[preprint,3p,12pt]{elsarticle}

\usepackage{amssymb}
\usepackage{lipsum,amsmath,algorithm,multirow}
\usepackage{amsfonts}
\usepackage{graphicx}
\usepackage{epstopdf}
\usepackage{algorithmic}
\usepackage{xcolor}
\usepackage{cleveref}
\usepackage{todonotes}
\usepackage{listofitems}
\usepackage{nicefrac}
\usepackage{xcolor,pgfplots}
\usepackage{rotating,booktabs,hvfloat,tikz}
\usepackage{url}
\usepackage{placeins}

\journal{Journal of Computational Physics}

\Crefname{figure}{Figure}{Figures}
\Crefname{table}{Table}{Tables}

\newcommand{\AH}[1]{{\color{black} #1}}
\newcommand{\AZ}[1]{{\color{black} #1}}

\pgfdeclareplotmark{myhalfcircle2}{
	\pgfpathcircle{\pgfpoint{0pt}{0pt}}{\pgfplotmarksize}
	\pgfusepathqstroke
	\pgfpathmoveto{\pgfpoint{\pgfplotmarksize}{0pt}}
	\pgfpatharc{0}{180}{\pgfplotmarksize}
	\pgfpathclose
	\pgfusepathqfill
}

\pgfdeclareplotmark{mymark1}{
	\pgfpathmoveto{\pgfqpoint{0pt}{\pgfplotmarksize}}%
	\pgfpathlineto{\pgfqpoint{.75\pgfplotmarksize}{0pt}}%
	\pgfpathlineto{\pgfqpoint{0pt}{-\pgfplotmarksize}}%
	\pgfpathlineto{\pgfqpoint{-.75\pgfplotmarksize}{0pt}}%
	\pgfpathclose
	\pgfusepathqstroke
	\pgfpathmoveto{\pgfqpoint{0pt}{-\pgfplotmarksize}}%
	\pgfpathlineto{\pgfqpoint{.75\pgfplotmarksize}{0pt}}%
	\pgfpathlineto{\pgfqpoint{-.75\pgfplotmarksize}{0pt}}%
	\pgfpathclose
	\pgfusepathqfill
}

\pgfdeclareplotmark{mymark2}{
	\pgfpathmoveto{\pgfqpoint{0pt}{\pgfplotmarksize}}%
	\pgfpathlineto{\pgfqpoint{.75\pgfplotmarksize}{0pt}}%
	\pgfpathlineto{\pgfqpoint{0pt}{-\pgfplotmarksize}}%
	\pgfpathlineto{\pgfqpoint{-.75\pgfplotmarksize}{0pt}}%
	\pgfpathclose
	\pgfusepathqstroke
	\pgfpathmoveto{\pgfqpoint{0pt}{\pgfplotmarksize}}%
	\pgfpathlineto{\pgfqpoint{.75\pgfplotmarksize}{0pt}}%
	\pgfpathlineto{\pgfqpoint{-.75\pgfplotmarksize}{0pt}}%
	\pgfpathclose
	\pgfusepathqfill
}

\begin{document}

\begin{frontmatter}



\title{Improving Pseudo-Time Stepping Convergence for CFD Simulations With Neural Networks}


\author[tub,comsol]{Anouk Zandbergen}
\ead{anouk.zandbergen@tu-berlin.de}
\affiliation[tub]{organization={Faculty V -- Mechanical Engineering and Transport Systems, Technische Universität Berlin },
            addressline={Straße des 17. Juni 135}, 
            city={Berlin},
            postcode={10623}, 
            country={Germany}}
        
\author[comsol]{Tycho van Noorden}
\ead{tycho.vannoorden@comsol.com}
\affiliation[comsol]{organization={COMSOL BV},
	addressline={Röntgenlaan 37}, 
	city={Zoetermeer},
	postcode={2719 DX}, 
	country={The Netherlands}}
        
\author[diam]{Alexander Heinlein}
\ead{a.heinlein@tudelft.nl}
\ead[url]{https://searhein.github.io/}
\affiliation[diam]{organization={Delft Institute of Applied Mathematics, Delft University of Technology},
	addressline={Mekelweg 4}, 
	city={Delft},
	postcode={2628 CD}, 
	country={The Netherlands}}

\begin{abstract}
Computational fluid dynamics (CFD) simulations of viscous fluids described by the Navier--Stokes equations are considered. Depending on the Reynolds number of the flow, the Navier--Stokes equations may exhibit a highly nonlinear behavior. The system of nonlinear equations resulting from the discretization of the Navier--Stokes equations can be solved using nonlinear iteration methods, such as Newton's method. However, fast quadratic convergence is typically only obtained in a local neighborhood of the solution, and for many configurations, the classical Newton iteration does not converge at all. In such cases, so-called globalization techniques may help to improve convergence.

In this paper, pseudo-transient continuation is employed in order to improve nonlinear convergence. The classical algorithm is enhanced by a neural network model that is trained to predict a local pseudo-time step. Generalization \AH{of the novel approach} is facilitated by predicting the local pseudo-time step separately on each element using only local information on a patch of adjacent elements as input. Numerical results for standard benchmark problems, including flow through a backward facing step geometry and Couette flow, show the performance of the machine learning-enhanced globalization approach; as the software for the simulations, the CFD module of COMSOL Multiphysics\textsuperscript{\tiny\textregistered} is employed.  
\end{abstract}



\begin{keyword}
Computational fluid dynamics \sep Newton's method \sep pseudo-time stepping \sep neural network \sep hybrid algorithm

\MSC[2020] 35Q30 \sep 65H10 \sep 65N12 \sep 65N22 \sep 68T07
\end{keyword}

\end{frontmatter}



\input{Mainmatter/Introduction}
\input{Mainmatter/NavierStokes}
\input{Mainmatter/PseudoTime}
\input{Mainmatter/NetworkForCFLpred/NetworkForCFLpred}
\input{Mainmatter/Results/Results}

\input{Mainmatter/Conclusion}


%
%
%

\appendix

\input{Mainmatter/Appendix}






\bibliographystyle{elsarticle-num} 
\bibliography{references}

\end{document}

%% file: Mainmatter/Introduction.tex

\section{Introduction}

Computational fluid dynamics (CFD) simulations are highly relevant for a wide range of applications, including the fields of aerospace, environmental, and biological engineering as well as weather predictions and medicine. Numerical simulations of Newtonian fluids\AH{, that is, fluid with a linear correlation of the local viscous stresses and strain rate,} involves the solution of the Navier--Stokes equations. Depending on the flow regime, the Navier--Stokes equations exhibit strong nonlinearities, which pose challenges to the numerical nonlinear solver employed.

In this work, we consider two-dimensional stationary CFD simulations using a mixed finite element method (FEM) as the discretization. The nonlinear system of equations resulting from discretizing the steady-state Navier--Stokes equations is then solved using Newton's method, that is, via solving a sequence of linearized problems. Whereas Newton's method provides quadratic convergence near the solution, the iteration may only converge slowly or even diverge farther away from the solution. In order to improve global convergence of the method, globalization techniques have been developed; see~\cite{MR2278446} for an overview of globalization techniques for the Navier--Stokes equations. In~\cite{MR2278446}, the authors categorize globalization techniques into the classes of \emph{backtracking methods} and \emph{trust region methods}; cf.~\cite{dennis_numerical_1996,eisenstat_globally_1994}. Other examples of methods to improve the robustness of the nonlinear solver include homotopy, continuation, pseudo transient continuation, mesh sequencing methods, or nonlinear preconditioning methods; cf.~\cite{allgower_continuation_1993,cliffe_numerical_2000,cresta_nonlinear_2007,hwang_parallel_2005,Kelley:1996:CAP,knoll_jacobian-free_2004,kubicek_computational_1983,liu_nonlinear_2022,watson_globally_1990,watson_algorithm_1997}.

The combination of scientific computing and machine learning, which is also denoted as scientific machine learning (SciML), is a new field~\cite{baker_workshop_2019}, which has gained a lot of attention over the past few years. One popular application area for SciML techniques are CFD simulations; see\AH{, for instance,} the review paper~\cite{brunton_machine_2020}. Successful examples for the use of SciML techniques in the context \AH{of CFD} are, for instance, discretization approaches using neural networks (NNs)~\cite{cai_physics-informed_2021,lagaris_artificial_1998,raissi_physics-informed_2019}, surrogate models~\cite{eichinger_surrogate_2022,guo_convolutional_2016,kim_deep_2019,kemna_reduced_2023}, and model discovery~\cite{beck_deep_2019,schmelzer_discovery_2020,brunton_discovering_2016}. In~\cite{wiewel_latent_2019}, the temporal evolution of dynamic systems is learned via long short-term memory (LSTM) NNs. Approaches for the enhancement of CFD simulations via machine learning include the generation of optimized meshes~\cite{huang_machine_2021}, the improvement of the resolution of the simulation results~\cite{kochkov_machine_2021,Marg:2022:ANN, Marg:2021:MultiGrid,margenberg_dnn-mg_2023}, or the detection \AH{of} troubled cells in simulation meshes due to Gibbs oscillations arising from the use of higher-order methods~\cite{RAY2019108845}. We note that, in~\cite{Marg:2022:ANN,Marg:2021:MultiGrid,margenberg_dnn-mg_2023,RAY2019108845}, the network inputs and output are specifically chosen from local patches of elements to achieve generalizability of the machine learning model; we will apply a similar approach in this work. 

The use of machine learning models for the enhancement of numerical solvers has also been investigated in several works. In the context of linear solvers, \cite{carlberg2016} and~\cite{pasetto_reduced_2017} have employed features generated by the proper orthogonal decomposition (POD) to develop a modified conjugate gradient (CG) solver or preconditioning techniques, respectively. \AH{The enhancement of the CG method by a nonlinear convolutional NN-based preconditioner has been investigate in~\cite{kaneda_deep_2022}.} In~\cite{heinlein_machine_2019,taghibakhshi_mg-gnn_2023}, NNs have been employed to improve domain decomposition methods, and in~\cite{Grebhahn:2017:PIM,antonietti_accelerating_2023,luz_learning_2020,moore_learning_2021,azulay_multigrid-augmented_2022}, components of multigrid methods are selected or generated using machine learning models. In the early works~\cite{Bhowmick:2006:AML,Bhowmick:2011:AAD,Bhowmick:2009:TLC,fuentes_statistical_2007,Motter:2015:LAS}, suitable numerical solution algorithms are being selected using machine learning techniques. Related to that, \cite{Holloway:2007:NNP} investigated the prediction of a \AH{suitable} combination of preconditioner and iterative method for solving a given sparse linear system.

Fewer works have considered the improvement of the nonlinear solver convergence using machine learning. In~\cite{luo_textpinmathcal_2023}, the nonlinear convergence of numerical solvers for nonlinear partial differential equations (PDE) is improved using POD features in a nonlinear preconditioning approach, and the proposed method is tested specifically for CFD problems. Moreover, in~\cite{Silva:2021:pmf}, dimensionless numbers and simulation properties are used as input to a random forest model to predict a relaxation parameter to accelerate a convergence of the Picard iteration for multiphase porous media flow. To the best of the authors' knowledge, there is currently no work that directly addresses the improvement of the nonlinear convergence of Newton's method for CFD simulations using NNs. The only exception is the master thesis of the first author~\cite{zandbergen_predicting_2022}, which laid the foundation for this paper.

In particular, we enhance the pseudo-transient continuation~\cite{Kelley:1996:CAP}, or pseudo-time stepping, approach for improving the robustness of the nonlinear convergence for the steady-state Navier-Stokes equations using an NN model. Similar to previous approaches~\cite{Marg:2022:ANN,Marg:2021:MultiGrid,margenberg_dnn-mg_2023,RAY2019108845}, the model is trained to predict the local pseudo-time step size for each mesh element using only local input features; this aims at ensuring generalizability of our approach. We test the performance and generalizability of our approach for different geometries and flow configurations using CFD simulations with the software package COMSOL Multiphysics\textsuperscript{\textregistered}; as a reference, we compare our approach against two default pseudo-time step control mechanisms implemented in the FEM software package COMSOL Multiphysics\textsuperscript{\textregistered}. As the machine learning framework, we employ Matlab, which enables us to access data from COMSOL via an available interface between the two packages. 

The paper is organized as follows. First, the Navier-Stokes equations and pseudo-time stepping are briefly recalled and discussed in~\cref{sec:NS} \AH{respectively}~\cref{sec:ptc}. In \cref{sec:NN} the NN for the local element pseudo-time step prediction is described. The numerical results are given in \cref{sec:results}. Lastly, the conclusions are presented in \cref{sec:conclusions}.

%% file: Mainmatter/NavierStokes.tex
\section{Navier--Stokes equations}
\label{sec:NS}

The Navier--Stokes equations are a system of partial differential equations (PDEs) that describes the flow of Newtonian fluids, that is, fluids with a linear correlation of the viscous stresses and the local strain rate. In particular, we consider incompressible fluids\AH{, that is, fluids with a constant density;} cf.~\cite{ECFD}. The stationary Navier-Stokes equations for incompressible flow read
\begin{align}
\begin{array}{rcl}
    \rho(\mathbf{u}\cdot\nabla)\mathbf{u} - \mu\nabla^2\mathbf{u} &=& -\nabla p + \mathbf{f}, \\
    \rho\nabla\cdot\mathbf{u} &=& 0,
    \end{array}\label{eq:NV}
\end{align}
describing the conservation of momentum and mass, respectively. In these equations, $\mathbf{u}$ is the flow velocity, $\mu$ is the dynamic viscosity of the fluid, $\rho$ is the density, which is assumed to be constant, $p$ is the pressure, and $\mathbf{f}$ is the body force acting on the flow \cite{ECFD}; in all our experiments, the body force is assumed to be zero. 

Due to the convective term $\rho(\mathbf{u}\cdot\nabla)\mathbf{u}$,~\cref{eq:NV} is nonlinear, and the numerical treatment of this nonlinearity is the subject of this research. Very strong nonlinearities are generally related to turbulent flow. However, already in the regime of laminar flow, the nonlinearities may become severe enough to cause nonlinear solvers to diverge; hence, we will focus on \AH{the laminar} case. Flow remains laminar as long as the Reynolds number is below a critical value in the order of $10^{5}$; see, e.g.,~\cite{incropera}. The Reynolds number $\text{Re}$ is given by
\begin{equation}
    \text{Re} = \frac{\rho UL}{\mu}.\label{eq:re}
\end{equation}
Here, $U$ and $L$ are the typical velocity and length scale. The Reynolds number describes the ratio between inertial and viscous forces; cf.~\cite{CMRM}. At low Reynolds numbers, the flow remains laminar as the viscous forces are able to damp out disturbances in the flow. At high Reynolds numbers, the inertial forces become more important resulting in nonlinear interactions to grow, which causes turbulence. 

In order to solve the stationary Navier-Stokes equations in~\cref{eq:NV} on a computational domain $\Omega \subset \mathbb{R}^2$ numerically, we consider a triangulation of $\Omega$ into linear triangles. Then, we discretize~\cref{eq:NV} using piecewise linear finite elements for both the velocity $\mathbf{u}$ and the pressure $p$ \AH{and stabilize} using Galerkin least-squares (GLS) streamline diffusion and Hughes–Mallet type shock-capturing \cite{HAUKE1994389,HUGHES1986305}. This yields a discrete but nonlinear system of equations
\begin{equation} \label{eq:discrete_NS}
\begin{aligned}
    N (\mathbf{u}_h) + B^\top p_h & = 0, \\
    B \mathbf{u}_h & = 0,
\end{aligned}
\end{equation}
where we denote the discrete velocity and pressure fields as $\mathbf{u}_h$ respectively $p_h$. The linear operators $B$ and $B^\top$ correspond to the divergence and gradient\AH{, respectively}, and the nonlinear operator $N$ corresponds to the convection-diffusion of the velocity field. 

By combining the fields $\mathbf{u}_h$ and $p_h$ into a single field $\mathbf{v}$, we can simply formulate the discrete nonlinear problem~\cref{eq:discrete_NS} as
$$
    F(\mathbf{v})=0.
$$
We solve this system using Newton's method, that is, by solving a sequence of linear systems of the form
$$
    -F'(\mathbf{v}^n) \Delta \mathbf{v}^n = F(\mathbf{v}^n),
$$
where $n$ is the iteration index of the Newton iteration \AH{and $F'(\mathbf{v}^n)$ is the Jacobian at the linearization point $\mathbf{v}^n$}. Close to the solution, Newton's method converges with quadratic rate. However, for an arbitrary initial guess, the method might also diverge. As we will see in our numerical results, this may happen for a relatively wide range of cases for the stationary Navier--Stokes equations. In order to make the nonlinear iteration more robust, globalization techniques can be employed; see, e.g.,~\cite{MR2278446} for an overview of globalization techniques for Newton's method for the Navier--Stokes equations. Here, we will consider the pseudo-time stepping approach~\cite{Kelley:1996:CAP}, which we will discuss in more detail in the next section.

%% file: Mainmatter/PseudoTime.tex
\section{Pseudo-time stepping}
\label{sec:ptc}

Solving stationary nonlinear equations \AH{using} Newton's method requires an initial guess close enough to the root. Globalization techniques are designed to give a result for a wider range of initial guesses, but can stagnate at local minima \cite{Kelley:1996:CAP}. Pseudo-time stepping, or pseudo-transient continuation, is an alternative to Newton's method for computing stationary solutions of time-dependent partial differential equations. \AH{The idea is to} obtain a solution more robustly by \AH{solving the time-dependent equation instead of the stationary one.} The idea behind the method is to frame the problem in a time-depending setting,
\begin{equation}
    \frac{\partial\mathbf{v}}{\partial t} = -F(\mathbf{v}),\label{eq:pseudotime}
\end{equation}
with $F(\mathbf{v}) = 0$ the system of nonlinear equations of which a solution must be determined. Now the algorithm for pseudo-time stepping can be described as the numerical integration with a variable time step method of the initial value problem 
\begin{equation}
    \frac{\partial\mathbf{v}}{\partial t}= -F(\mathbf{v}),\ \ \ \mathbf{v}(0) = \mathbf{v}_0; 
    \label{eq:pseudotime2}
\end{equation}
\AH{cf.~\cite{Kelley:1996:CAP}.} The idea is that this time integration converges to a steady state for a wider range of initial values $\mathbf{v}_0$ than \AH{the} standard Newton method. In addition, the method tries to increase the time step as $F(\mathbf{v})$ approaches 0, such that when the iterations get closer to a steady state, the convergence becomes (close to) quadratic.

The iterations of the pseudo-time stepping algorithm are given by
\begin{equation}
    \mathbf{v}^{n+1} = \mathbf{v}^n - \left((\Delta t^n)^{-1}I + F'(\mathbf{v}^n)\right)^{-1}F(\mathbf{v}^n),\label{eq:iteration}
\end{equation}
with $F'(\mathbf{v}^n)$ the Jacobian. To understand how \eqref{eq:iteration} is obtained from \eqref{eq:pseudotime2} consider an Euler backward step starting from $\mathbf{v}^n$ for equation \eqref{eq:pseudotime2}:
\begin{equation}
    \mathbf{z}^{n+1} = \mathbf{v}^n - \Delta t^n F(\mathbf{z}^{n+1}).
\end{equation}
Note that $\mathbf{z}^{n+1}$ is a root of
\begin{equation}
    G(\xi): = \xi + \Delta t^n F(\xi)-\mathbf{v}^n,
\end{equation}
as a function of $\xi$. The first Newton iterate for solving $G(\xi)=0$ with initial guess $\xi_0 = \mathbf{v}^n$ gives
\begin{align*}
    \xi_1 &= \mathbf{v}^n - \left(I+\Delta t^n F'(\mathbf{v}^n)\right)^{-1}(\mathbf{v}^n+\Delta t^n F(\mathbf{v}^n)-\mathbf{v}^n), \\
    &= \mathbf{v}^n - \left((\Delta t^n)^{-1}+F'(\mathbf{v}^n)\right)^{-1}F(\mathbf{v}^n).
\end{align*}
From this corrector iteration, we have obtained the formula in \eqref{eq:iteration}. The \AH{pseudo code of the method is the given by}~\Cref{alg:pts}.

\begin{algorithm}
	\caption{Pseudo-time stepping \AH{algorithm from}~\cite{Kelley:1996:CAP}} \label{alg:pts}
	\begin{algorithmic}[1]
	\STATE Set $\mathbf{v}=\mathbf{v}_0$ and $\Delta t = \Delta t_0$.
	\WHILE {$\|F(\mathbf{v})\|$ is too large}
		\STATE Solve $(\Delta t^{-1}I+F'(\mathbf{v}))\mathbf{s} = -F(\mathbf{v})$
		\STATE Set $\mathbf{v}=\mathbf{v}+\mathbf{s}$
		\STATE Evaluate $F(\mathbf{v})$
		\STATE Update $\Delta t$
		\ENDWHILE
	\end{algorithmic} 
\end{algorithm}

\subsection{Applying pseudo-time stepping to the weak form of the Navier-Stokes equations}
The weak form of the Navier-Stokes equations is given by: find velocity and pressure fields $(\mathbf{u},p)\in X$ such that
\begin{align}
\begin{array}{rcl}
    \displaystyle{\int_\Omega} (\rho\partial_t\mathbf{u}+\rho(\mathbf{u}\cdot\nabla)\mathbf{u})\cdot\phi - 
    (\mu\nabla\mathbf{u}\cdot\nabla\phi-p\nabla\cdot\phi)&=& \displaystyle{\int_\Omega}\mathbf{f}\cdot\phi+\int_{\Gamma_N} \mathbf{h}\cdot \phi, \\
    \displaystyle{\int_\Omega}\rho\nabla\cdot\mathbf{u}\psi &=& 0,
    \end{array}\label{eq:NV2}  
\end{align}
for all test functions $(\phi,\psi)\in Y$, with $X$ and $Y$ suitable function spaces, e.g., $X = Y = H_0^1(\Omega)^2 \times L^2(\Omega)$. Here $\mathbf{f}$ is a forcing term, and $\mathbf{h}$ encodes the boundary conditions.

To apply pseudo-time stepping to this weak form, we first discretize in space using finite elements by replacing the function spaces $X$ and $Y$ by finite dimensional subspaces $X_h\subset X$ and $Y_h\subset Y$, consisting of continuous functions that are piecewise linear on a finite element mesh of the domain $\Omega$ with maximum mesh size $h$. Then, in order to apply pseudo-time stepping to the stationary form of the spatially discretized system, the same conceptual idea as described in the previous section is used, which can be summarized by first writing the time-discretized equations for one Euler backward step, and then linearizing this system for the new time step around the current time step. The Euler backward formula with time step $\Delta t^n$ for the discretized weak form of \eqref{eq:NV2} gives the following (discretized) nonlinear problem: find $(\mathbf{u}_h,p_h)\in X_h$ such that 
\begin{align}
\begin{array}{r}
    \displaystyle{\sum_{e\in E}}\Bigg(\displaystyle{\int_{\Omega^e}} (\rho \frac{\mathbf{u}_h^{n+1}-\mathbf{u}_h^{n}}{\Delta t^n}+\rho(\mathbf{u}_h^{n+1}\cdot\nabla)\mathbf{u}_h^{n+1})\cdot\phi_j - 
    (\mu\nabla\mathbf{u}_h^{n+1}\cdot\nabla\phi_j - p_h^{n+1}\nabla\cdot\phi_j) + \\
    -\displaystyle{\int_{\Omega^e}}\mathbf{f}\cdot\phi_j-\displaystyle{\int_{\Gamma^e_N}} \mathbf{h}\cdot \phi_j\Bigg)=0, \\
    \displaystyle{\sum_{e\in E}}\Bigg(\displaystyle{\int_{\Omega^e}}\rho\nabla\cdot\mathbf{u}_h^{n+1}\psi_j\Bigg) = 0,
    \end{array}
\label{eq:NVEB}   
\end{align}
for all test functions $(\phi_j,\psi_j) \in Y_h$, where the summation is over the set of \AH{all} mesh elements $E=\{e_1,e_2,\ldots, e_m\}$, with the number of elements \AH{$m$}.

Linearizing this system around $(\mathbf{u}_h^n,p_h^n)$ results in
\begin{align}
\begin{array}{r}
\displaystyle{\sum_{e\in E}}\Bigg(
    \displaystyle{\int_{\Omega^e}} (\rho \frac{\mathbf{u}^{n+1}_h-\mathbf{u}^{n}_h}{\Delta t^n}+\rho(\mathbf{u}^{n}_h\cdot\nabla)\mathbf{u}^{n+1}_h+\rho(\mathbf{u}^{n+1}_h\cdot\nabla)\mathbf{u}^{n}_h-\rho(\mathbf{u}^{n}_h\cdot\nabla)\mathbf{u}^{n}_h)\phi_j+\\ - \displaystyle{\int_{\Omega^e}}(\mu\nabla\mathbf{u}^{n+1}_h\cdot\nabla\phi_j-p^{n+1}_h\nabla\cdot\phi_j)- \displaystyle{\int_{\Omega^e}}\mathbf{f}\cdot\phi_j-\displaystyle{\int_{\Gamma^e_N}} \mathbf{h}\cdot \phi_j\bigg)=0, \\
    \displaystyle{\sum_{e\in E}}\Bigg(\displaystyle{\int_{\Omega^e}}\rho\nabla\cdot\mathbf{u}^{n+1}_h\psi_j \bigg)=0.
    \end{array}\label{eq:NVPT}   
\end{align}
Solving this linear system for $(\mathbf{u}^{n+1}_h,p^{n+1}_h)$ constitutes one pseudo-time step. We write~\cref{eq:NVEB} in matrix format
\begin{align}
\frac{1}{\Delta t^n}M(\mathbf{v}^{n+1}-\mathbf{v}^n)+F(\mathbf{v}^{n+1})=0,    
\end{align}
where the coefficients of $(\mathbf{u}^n_h, p^n_h)$ w.r.t.\ the chosen basis of $X_h$ are again combined into a single vector $\mathbf{v}^n$, and where $M = \begin{bmatrix}
    M_{\mathbf{u}} & \mathbf{0} \\ \mathbf{0} & \mathbf{0}
\end{bmatrix}$ is a singular mass matrix with the velocity mass matrix $M_{\mathbf{u}}$. Then, we can also write the solution of~\cref{eq:NVPT} in the more compact form
\begin{align}
\mathbf{v}^{n+1}=\mathbf{v}^n-(\frac{1}{\Delta t^n}M+F'(\mathbf{v}^n))^{-1}F(\mathbf{v}^n). \label{eq:LNVPT}
\end{align}

\subsection{Choice of the pseudo-time step} \label{sec:local_cfl}

The idea \AH{is that} the pseudo-time step $\Delta t^n$ \AH{is initially small}, \AH{such that its influence in}~\eqref{eq:iteration} is large. \AH{W}hen the \AH{iteration} \AH{approaches} convergence, \AH{the magnitude of the pseudo-time steps becomes} small, \AH{such that}~\eqref{eq:iteration} performs more like \AH{the} standard Newton iteration. This has the benefit that\AH{, in the beginning,} the method \AH{depends less} on the initial guess\AH{, and in the end, it may converge} quadratically. There are several possibilities to define a pseudo-time step that fits this property.

In addition \AH{to that,} it is possible to use a pseudo-time step that uses local, i.e.\ mesh element dependent, information. This amounts to replacing the global pseudo-time step $\Delta t^n$ in \eqref{eq:NVPT} by a local, mesh element dependent, pseudo-time step $\Delta t^n_e$. In this case the mass matrix $M$ in \eqref{eq:LNVPT} will depend on the vector of local pseudo-time steps $\overline{\Delta t}^n=(\Delta t^n_{e_1},\ldots,\Delta t^n_{e_m})$, so that \eqref{eq:LNVPT} becomes
\begin{align}
\mathbf{v}^{n+1}=\mathbf{v}^n-(M(\overline{\Delta t}^n)+F'(\mathbf{v}^n))^{-1}F(\mathbf{v}^n). \label{eq:LNVPTL}
\end{align}

In COMSOL, for instance, the following local pseudo-time step is used\AH{:}
\begin{equation}
    \Delta t^n_e = \text{CFL}(n)\frac{h_e}{\|\mathbf{u}^n_h\|_e},\label{eq:deltaT}
\end{equation}
with $h_e$ the size of mesh element $e$, $\|\mathbf{u}^n_h\|_e$ the Euclidean norm of the velocity field $\mathbf{u}^n_h$ in the mesh element $e$, and CFL$(n)$ a global Courant--Friedrichs--Lewy (CFL) number~\cite{CMRM}. This global CFL number \AH{may depend} on the iteration count $n$, \AH{for instance,} 
\begin{equation}\label{eq:CFLloc} 
    \text{CFL}_{iter}(n) 
    = 
    \begin{cases} 1.3^{\min(n,9)} & 1\leq n \leq 20, \\ 
    	1.3^9+9\cdot 1.3^{\min(n-20,9)} & 20 < n \leq 40, \\ 
	    1.3^9+9\cdot 1.3^9+ 90\cdot1.3^{\min(n-40,9)} & n > 40\AH{;}
    \end{cases}
\end{equation}
\AZ{see \cite[p.~1108]{CMRM}}. \AH{Alternatively,} it can be given by a controller based on the nonlinear error estimate $e_n$ for iteration $n$, the given target error estimate ``tol'', and \AH{control} parameters $k_P$, $k_I$, and $k_D$, which are positive constants:
\begin{equation}\label{eq:CFL2}
    \text{CFL}_e(n)=\left(\frac{e_{n-2}}{e_{n-1}}\right)^{k_P} \left(\frac{\text{tol}}{e_{n-1}}\right)^{k_I}\left(\frac{e_{n-2}/e_{n-1}}{e_{n-3}/e_{n-2}}\right)^{k_D}\text{CFL}_e(n-1),
\end{equation}
where $e_n$ is an estimate for the error in the $n$th nonlinear iteration; see~\cite[p.~1515]{CMRM}. Both these options are available for controlling the global CFL number in COMSOL. 

\begin{figure}[t]
	\centering
	\includegraphics[width=0.4\textwidth]{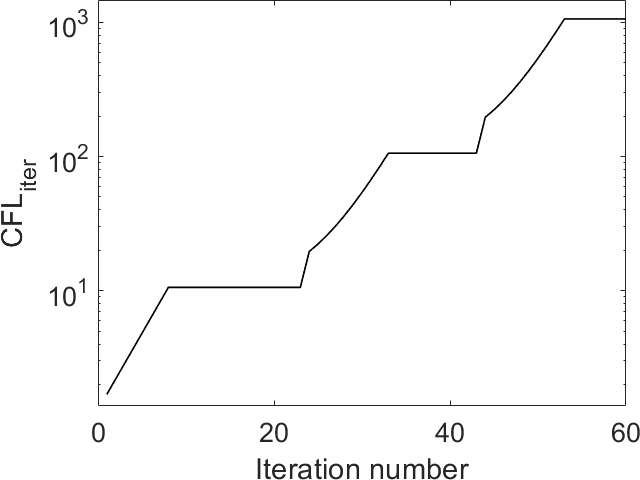}
	\caption{A plot of CFL$_{iter}$ \AH{depending} on the nonlinear iteration count\AZ{, as \AH{defined} in~\cref{eq:CFLloc}}.}
	\label{fig:cflnum}
\end{figure}

For both cases, in the initial iterations of pseudo-time stepping, the method starts with a small global CFL number and gradually increases it as the solution reaches convergence~\cite{CMRM}. In~\cref{eq:CFLloc} this is because the CFL number increases each iteration, also shown in~\Cref{fig:cflnum}. For~\cref{eq:CFL2}, when the solution is near convergence, the error estimates $e_{n}$ will be smaller and thus the CFL number increases. 

In this paper we investigate replacing these two pseudo-time step control algorithms by the use of an NN to control the local, mesh element dependent, pseudo-time step $\Delta t^n_e$. In order to make the method easily generalizable to different meshes, the NN is provided with local information from \AH{the corresponding} mesh element and its adjacent mesh elements. One could think of local information about the solution, residuals, mesh and the cell Reynolds number. For example, elements with high local residuals could then have smaller local pseudo-time steps compared to elements with low residuals, which could accelerate convergence in areas which already have low residuals. So, treating each mesh element individually using local information can accelerate convergence for the whole simulation.

%% file: Mainmatter/NetworkForCFLpred/NetworkForCFLpred.tex
\section{Neural network for the local pseudo-time step prediction}
\label{sec:NN}
\AH{In order to} predict an element-wise local pseudo-time step for the pseudo-time stepping that yields fast and robust Newton convergence, we employ a data-based approach based on machine learning techniques; more specifically, we use artificial neural networks (ANNs). The most basic form of an ANN is a multilayer perceptron (MLP) and is given as a composed function
$$
\begin{aligned}
    \mathcal{NN} \left( x \right) 
    =
    W_{n+1} \ f_{n} \circ \dots \circ f_{1} \left( x \right),
\end{aligned}
$$
with
\begin{equation} \label{eq:nn}
	f_i \left( x \right) = \sigma \left( W_i x + b_i \right).
\end{equation}
Here, \AH{the} $W_i \in \mathbb{R}^{n_{i} \times n_{i-1}}$ and $b_i \in \mathbb{R}^{n_{i}}$ are the so-called weight matrices and bias vectors, which represent the linear parts of the NN function $\mathcal{NN}$; see~\Cref{fig:nn1} for an exemplary network architecture of an MLP. \AH{G}iven a specific network architecture, that is, the dimensions $\left(n_i\right)_i$, the MLP is determined by the coefficients of the weight matrices and bias vectors, also denoted as the network parameters. The MLP becomes nonlinear due to composition with the activation function $\sigma$. Here, we will employ the rectified linear unit (ReLU) function $\sigma \left( x \right) = \max \left( 0 , x \right)$. MLPs, and NNs in general, are well-suited for approximating nonlinear functions; see, e.g.,~\cite{cybenko_approximation_1989}. In order fit an NN to input data $X = \left\lbrace x_i \right\rbrace$ and corresponding output data $Y = \left\lbrace y_i \right\rbrace$, a loss function $\mathcal{L}$ is minimized with respect to the network parameters:
$$
    \arg\min_{W_i,b_i} \mathcal{L} \left( \mathcal{NN}\left( X \right),Y \right)
$$
The loss function penalizes deviation of the network function from the reference output data $Y$. In order to optimize the NN parameters (a.k.a.~network training), we employ a stochastic gradient descent (SGD) method using adaptive moment estimation (Adam)~\cite{kingma_adam_2017}\AH{; the gradients are computed using the backpropagation algorithm~\cite{kelley1960gradient}.} For a more detailed introduction to deep learning and NNs; see\AH{, e.g.,}~\cite{goodfellow_deep_2016}.

In this section, we will discuss \AH{how} the local pseudo-time steps used as reference data in the root mean squared error (RMSE) loss function \AH{are computed as well as} the composition of our training data set \AH{and} the network architecture. 

\input{Mainmatter/NetworkForCFLpred/LOSS}

\input{Mainmatter/NetworkForCFLpred/DataSets}

\input{Mainmatter/NetworkForCFLpred/TrainingData}

\input{Mainmatter/NetworkForCFLpred/NetworkStructure}

%% file: Mainmatter/NetworkForCFLpred/LOSS.tex
\subsection{Loss and optimal local pseudo-time step}\label{sec:loss}

\begin{figure}[t]
	\centering
	\includegraphics[width=0.9\textwidth,trim=0mm 0mm 0mm 4mm,clip]{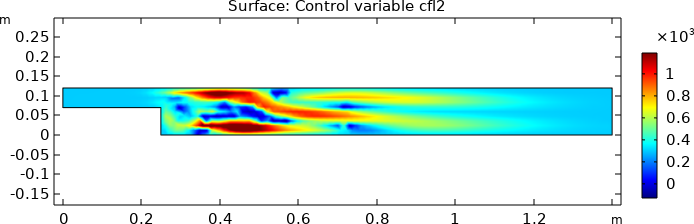}
	\caption{The local pseudo-time step \AH{optimized by SNOPT; the difference between} the obtained solution and the converged solution was minimized \AH{up} to a value of $3.1{\cdot}10^{-5}$.}
	\label{fig:CFLCOMSOL}
\end{figure}

Our goal is to predict local pseudo-time steps that enhance the pseudo-time stepping method in accelerating convergence of a given flow problem. In particular, we try to predict the vector of local pseudo-time steps $\overline{\Delta t}_{opt}$ such that the next iterate is as close as possible to the solution of the nonlinear problem. We will denote this as the \textit{the optimal local pseudo-time step} for an iterate $\mathbf{v}^n$, $\overline{\Delta t}_{opt}(\mathbf{v}^n)$, defined by
\begin{equation}
   \overline{\Delta t}_{opt}(\mathbf{v}^n) := \underset{\overline{\Delta t}}{\arg\min}\, || \mathbf{v}^n - \left(M(\overline{\Delta t}) +F'(\mathbf{v}^n)\right)^{-1}F(\mathbf{v}^n)-\mathbf{v}^* ||,\label{eq:optcfl}
\end{equation}
where $\mathbf{v}^*$ satisfies $F(\mathbf{v}^*)=0$. The norm that is used here is the $L^2$ norm of the reconstructed velocity components.

Since the optimal local pseudo-time step cannot be computed explicitly, we compute an approximation via an optimization procedure, using the sparse nonlinear optimizer (SNOPT) software package~\cite{philip2015user,Gill:2002:SNOPT}. Th\AH{is} optimization is performed within COMSOL\AH{; see~\cite{zandbergen_predicting_2022} for more details.} For the example of a back-step geometry with laminar flow, we obtain the optimal local pseudo-time step as shown in~\Cref{fig:CFLCOMSOL} by using this two step procedure. 

Once we have computed the local optimal pseudo-time step for \AH{all elements in the data set}, we optimize the squared error of the network prediction against \AH{these} reference value\AH{s}. The whole data set \AH{consists} of a large number of optimal local pseudo-time steps $\Delta t_{e_i,opt}(\mathbf{v}^n)$ computed from various problem configurations, nonlinear iterates $\mathbf{v}^n$, and elements $e_i$ computed using the optimization in COMSOL. To define the complete loss, let $\Delta t^{(i)}_{pred}$ and $\Delta t_{opt}^{(i)}$ be the network prediction and optimized pseudo-time steps, respectively, for the $i$th data point in a data set with $n$ data points. The full RMSE loss, is then given as:
\begin{equation}\label{eq:rmseloss}
    \text{loss} = \sqrt{\frac{1}{n}\sum_{i = 1}^n\left( \Delta t_{pred}^{(i)} - \Delta t_{opt}^{(i)}  \right)^2}.
\end{equation}

%% file: Mainmatter/NetworkForCFLpred/DataSets.tex
\subsection{Neural network input} \label{sec:DataSets}

\begin{figure}[t]
	\centering
	\input{Mainmatter/NetworkForCFLpred/patch}
	\caption{Data point which contains information of a patch of four elements $e_{j,1}$ to $e_{j,4}$ with their vertices $v_{j,1}$ to $v_{j,6}$\AH{; see~\Cref{tab:input} for all input features for the NN model on each element patch}.
	\label{fig:patch}}
\end{figure}
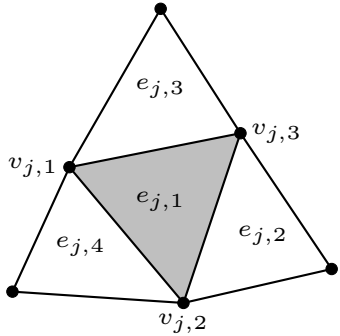

\begin{table}[t]
	\centering 
	\begin{tabular}{|l|l|l|}
		\hline
		\textbf{variable}    & \textbf{unit}    & \textbf{location}              \\ \hline\hline
		element edge length  & m                & -                              \\ \hline
		$u$                  & \nicefrac{m}{s}              & $v_{j,1}$, $v_{j,2}$, $v_{j,3}$ \\ \hline
		$v$                  & \nicefrac{m}{s}               & $v_{j,1}$, $v_{j,2}$, $v_{j,3}$ \\ \hline
		$p$                  & Pa               & $v_{j,1}$, $v_{j,2}$, $v_{j,3}$ \\ \hline
		$R_u$                &                  & $v_{j,1}$, $v_{j,2}$, $v_{j,3}$ \\ \hline
		$R_v$                &                  & $v_{j,1}$, $v_{j,2}$, $v_{j,3}$ \\ \hline
		$R_p$                &                  & $v_{j,1}$, $v_{j,2}$, $v_{j,3}$ \\ \hline
		$r_u$                & \nicefrac{N}{m$^3$} & $v_{j,1}$, $v_{j,2}$, $v_{j,3}$ \\ \hline
		$r_v$                & \nicefrac{N}{m$^3$}          & $v_{j,1}$, $v_{j,2}$, $v_{j,3}$ \\ \hline
		$r_p$                & \nicefrac{kg}{(m$^3\cdot$s)} & $v_{j,1}$, $v_{j,2}$, $v_{j,3}$  \\ \hline
		cell Reynolds number &                  & $e_{j,1}$                      \\ \hline
	\end{tabular}\caption{Variables used for input for the NN with their given unit and location on one element in the patch, in this case the central element. One data point contains information of four elements, with a non-structured order of the adjacent elements. The two different types of residuals $r_\alpha$ and $R_\alpha$, $\alpha= u,\,v,\,p$ are given in~\cref{eq:R,eq:r}.
	\label{tab:input}}
\end{table}

NNs are universal function approximators and can theoretically approximate any continuous nonlinear function up to arbitrary precision; see, for instance,~\cite{cybenko_approximation_1989,hornik_multilayer_1989,hornik_approximation_1991}. Therefore, we expect that it is possible to predict the optimal local pseudo-time step, given sufficient input data. On the other hand, we want to limit the complexity of the input data and train a model with strong generalization properties. Therefore, similar to other approaches for NN-enhanced simulations~\cite{Marg:2022:ANN,Marg:2021:MultiGrid,RAY2019108845}, we employ local input data to predict the local optimal pseudo-time step. This means that we do not include any information about the specific global boundary value problem but only data that is given on a local patch of elements around the element of interest\AH{, that is}, where the local optimal pseudo-time step is to be predicted. As a result, the trained NN model is applicable to any considered laminar flow problem. An additional benefit of using only local data is that the dimensionality and complexity of the data is being limited. \AH{Note that,} as also noted in~\cite{Marg:2022:ANN}, \AH{due to the local data approach,} a single simulation already generates a rich data set: each simulation yields data on a large set of elements over a number of Newton iterations. However, to generate a training data set with sufficient variation to allow for generalization of the model, the training data may be sampled from different simulations with varying flow conditions, induced by different boundary conditions and geometries as well as mesh refinement levels. 

Let us now discuss the specific data included in the model input. As mentioned before, we employ data from a patch of elements
$$
	P_j = \bigcup_{i=1}^4 e_{j,i}.
$$ 
In particular, the patch consists of the element for which we want to compute the local optimal pseudo-time step, $e_{j,1}$, as well as the \AH{(up to)} three adjacent elements $e_{j,2}$, $e_{j,3}$, and $e_{j,4}$; cf.~\cref{fig:patch}. Specifically, we denote two elements as adjacent if they are connected via an element edge. On each patch, we sample the data listed in~\cref{tab:input} with their respective location. The nodal information is explicitly available from the mixed finite element discretization used in our simulations; cf.~the description of our simulation setup in~\cref{sec:NS}.

Let us briefly motivate our choices. Firstly, we include the element size, which is implied by the element edge lengths, and $\mathbf{u}$ since they are also employed in the computation of the pseudo-time step approach in~\cref{eq:deltaT}. We complement this $\mathbf{u}$ by $p$ to provide a complete description of the finite element solution on the local patch. Moreover, the local information describing the convergence of the Newton iteration is supplied in the form of two different types of residuals, i.e., the residual \AH{of} the discretized system of nonlinear equations
\begin{equation}
F(\mathbf{v})
= 
\begin{pmatrix} R_{u} \\ R_{v} \\ R_p \end{pmatrix} 
= 
\begin{pmatrix}
N(u,v) +B^\top p \\
B(u,v)
  \end{pmatrix}
  = \begin{pmatrix}
N(\mathbf{u}) +B^\top p \\
B(\mathbf{u})
  \end{pmatrix}
\label{eq:R}
\end{equation}
and the residual obtained by substituting the approximate solution into the system of PDEs~\cref{eq:NV}:
\begin{equation}
\begin{aligned}
r_u&=\rho(u\partial_xu+v\partial_y u)-\mu(\partial_x^2 u+\partial_y^2 u)+\partial_x p+f_x, \\
r_v&=\rho(u\partial_xv+v\partial_y v)-\mu(\partial_x^2 v+\partial_y^2 v)+\partial_y p+f_y, \label{eq:r}\\
r_p&=\rho(\partial_xu+\partial_y v). \end{aligned}
\end{equation}

\AH{The} local Reynolds number\AH{, as given} in equation~\cref{eq:re}, provides additional information about the local nonlinearity of the flow; it is also used in papers \cite{Marg:2022:ANN,Marg:2021:MultiGrid,Silva:2021:pmf} as input for the machine learning model. The input is normalized to have mean 0 and standard deviation by computing the z-score \AZ{defined as
$$
    z(\sigma) = \frac{x-\mu}{\sigma},
$$
with $x$ the evaluated input data, $\mu$ the mean of the input, and $\sigma$ the standard deviation of the input \cite{kreyszig2010advanced}}. The the mean and standard deviation for centering and scaling the data are saved and later on used to normalize input data from other simulations.

Note that, for practical reasons\AH{,} in our implementation in COMSOL, we collect all input information, as listed in~\cref{tab:input}, element by element from each patch. This means that vertex- and edge-based information is duplicated where the elements in the patch touch. As a result, we obtain a total of $n_0 = 124$ input features for our network model.
In practice, one may want to omit the redundant data, however, this is \AH{currently} not straightforward based on the data structures in our interface from Matlab to COMSOL since the elements in each patch cannot be easily retrieved in a consistent ordering. \AH{Moreover}, the redundant information implicitly encodes connectivity information, which could alternatively be encoded by a consistent ordering of the elements within the patch. Our numerical results in~\cref{sec:results} indicate that this handling of the input data does not prevent our model to learn from the data. However, in future work, we plan to investigate if using unique input data with a consistent ordering may have a positive effect on our model. 

Finally, elements which are adjacent to the boundary of the computational domain $\Omega$, may only have one or two neighboring elements, resulting in patches consisting of two or three elements only. In order to encode this in the data, we set any data of the vertices and edges outside the domain to zero; other approaches for encoding this information are possible and will, again, be considered in future work. 

%% file: Mainmatter/NetworkForCFLpred/patch.tex
\begin{tikzpicture}[scale = 1.5, transform shape]
\draw[black,thick,fill=gray!50] (3,0) -- (3.5,1.5) -- (2,1.2) -- (3,0) -- cycle;
\draw[black,thick] (3,0) -- (1.5,0.1) -- (2,1.2);
\draw[black,thick] (2,1.2) -- (2.8,2.6) -- (3.5,1.5);
\draw[black,thick] (3,0) -- (4.3,0.3) -- (3.5,1.5);
\draw[black,fill=black] (3,0) circle(0.05cm) 
node[node font = \tiny,below]{$v_{j,2}$};
\draw[black,fill=black] (3.5,1.5) circle(0.05cm)
node[node font = \tiny,anchor = west]{$v_{j,3}$};
\draw[black,fill=black] (2,1.2) circle(0.05cm)
node[node font = \tiny,anchor = east]{$v_{j,1}$};
\draw[black,fill=black] (1.5,0.1) circle(0.05cm);
\draw[black,fill=black] (2.8,2.6) circle(0.05cm);
\draw[black,fill=black] (4.3,0.3) circle(0.05cm);
\draw (2.8,0.9) node[node font = \tiny]{$e_{j,1}$};
\draw (2.1,0.5) node[node font = \tiny]{$e_{j,4}$};
\draw (3.7,0.6) node[node font = \tiny]{$e_{j,2}$};
\draw (2.8,1.9) node[node font = \tiny]{$e_{j,3}$};
\end{tikzpicture}

%% file: Mainmatter/NetworkForCFLpred/TrainingData.tex
\subsection{Generation of training data}
\label{sec:traindata}

\begin{table}[t]
	\centering \small
	\begin{tabular}{|l|r|r|r|r|}
		\hline
		\textbf{back-step} &                     \textbf{dim. inflow} &         \textbf{dim. outflow} &              \textbf{range inflow} &                 \textbf{range maximum} \\
		\textbf{geometry}  &                      \textbf{tunnel (m)} &           \textbf{tunnel (m)} &            \textbf{velocity (\nicefrac{m}{s})} &              \textbf{element size (m)} \\ \hline\hline
		\textbf{B1}        & 0.05\phantom{0} $\times$ 0.25\phantom{0} & 0.12 $\times$ 1.15\phantom{0} &                     0.001 -- 0.015 & 0.0106\phantom{0} -- 0.0256\phantom{0} \\ \hline
		\textbf{B1S}        &                     0.005 $\times$ 0.025 &          0.012 $\times$ 0.115 & 0.01\phantom{0} -- 0.15\phantom{0} &                     0.00106 -- 0.00256 \\ \hline
		\textbf{B2}        & 0.08\phantom{0} $\times$ 0.25\phantom{0} & 0.22 $\times$ 1.15\phantom{0} &           0.001 -- 0.01\phantom{5} & 0.0126\phantom{0} -- 0.0206\phantom{0} \\ \hline
		\textbf{B2S}        &                     0.008 $\times$ 0.025 &          0.022 $\times$ 0.115 & 0.01\phantom{0} -- 0.1\phantom{05} &                     0.00126 -- 0.00206 \\ \hline
	\end{tabular}
	\caption{\AH{Dimensions of the back-step geometries considered for the training and test data; see~\Cref{fig:backstep} for the base geometry.}
	\label{tab:backstepgeom}}
\end{table}

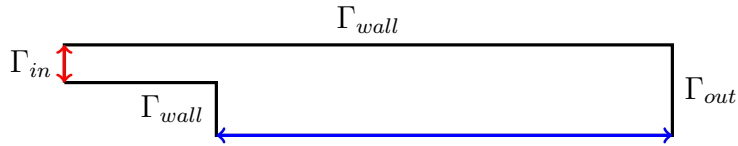
\begin{figure}
	\centering
	\input{Mainmatter/NetworkForCFLpred/backstep}
	\caption{
		\AH{Back-step geometries employed: B1, B1S, B2, and and B2S vary in the width of the inlet $\Gamma_{in}$ (red), where a parabolic velocity profile is prescribed, and the length of the blue wall segment; cf.~\Cref{tab:backstepgeom}. At the outlet $\Gamma_{out}$, we prescribe a constant pressure boundary condition, and at all other parts of the boundary, which are denoted as $\Gamma_{wall}$, a no-slip boundary condition is enforced. For an exemplary flow field, see~\Cref{fig:b1}.}
	\label{fig:backstep}}
\end{figure}

\begin{table}[t!]
	\centering \small
	\begin{tabular}{|l|r|r|r|r|}
		\hline
		\textbf{back-step}              				 &       \textbf{number of} &   \textbf{max.\ element} & \textbf{inflow velocity} & \multirow{2}{*}{\textbf{\# iterations}} \\
		\textbf{geometry}                                                 &        \textbf{elements} &        \textbf{size (m)} &           \textbf{(\nicefrac{m}{s})} &                                      \\ \hline\hline
		\multirow{12}{*}{\textbf{B1}}                    &  \multirow{3}{*}{6\,516} &  \multirow{3}{*}{0.0156} &    {\footnotesize 0.001} &                {\footnotesize 1, 10} \\
		                                                 &                          &                          &    {\footnotesize 0.005} &                   {\footnotesize 10} \\
		                                                 &                          &                          &     {\footnotesize 0.01} &                   {\footnotesize 10} \\ \cline{2-5}
		                                                 &  \multirow{3}{*}{3\,908} &  \multirow{3}{*}{0.0206} &    {\footnotesize 0.003} &                {\footnotesize 1, 10} \\
		                                                 &                          &                          &    {\footnotesize 0.006} &                   {\footnotesize 10} \\
		                                                 &                          &                          &    {\footnotesize 0.009} &                   {\footnotesize 10} \\ \cline{2-5}
		                                                 &  \multirow{4}{*}{3\,558} &  \multirow{4}{*}{0.0266} &    {\footnotesize 0.002} &                {\footnotesize 1, 10} \\
		                                                 &                          &                          &    {\footnotesize 0.007} &                   {\footnotesize 10} \\
		                                                 &                          &                          &    {\footnotesize 0.008} &                 {\footnotesize 2, \phantom{0}4} \\
		                                                 &                          &                          &    {\footnotesize 0.015} &                   {\footnotesize 10} \\ \cline{2-5}
		                                                 & \multirow{2}{*}{13\,828} &  \multirow{3}{*}{0.0106} &    {\footnotesize 0.004} &                {\footnotesize 2, 10} \\
		                                                 &                          &                          &    {\footnotesize 0.008} &                {\footnotesize 2, 10} \\ \hline
		\multirow{12}{1.8cm}{\textbf{B1S} \\ (scaled B1)} &  \multirow{3}{*}{6\,516} & \multirow{3}{*}{0.00156} &     {\footnotesize 0.01} &                {\footnotesize 1, 10} \\
		                                                 &                          &                          &     {\footnotesize 0.05} &                   {\footnotesize 10} \\
		                                                 &                          &                          &      {\footnotesize 0.1} &                   {\footnotesize 10} \\ \cline{2-5}
		                                                 &  \multirow{3}{*}{3\,908} & \multirow{3}{*}{0.00206} &     {\footnotesize 0.03} &                {\footnotesize 1, 10} \\
		                                                 &                          &                          &     {\footnotesize 0.06} &                   {\footnotesize 10} \\
		                                                 &                          &                          &     {\footnotesize 0.09} &                   {\footnotesize 10} \\ \cline{2-5}
		                                                 &  \multirow{4}{*}{3\,558} & \multirow{4}{*}{0.00266} &     {\footnotesize 0.02} &                {\footnotesize 1, 10} \\
		                                                 &                          &                          &     {\footnotesize 0.07} &                   {\footnotesize 10} \\
		                                                 &                          &                          &     {\footnotesize 0.08} &                 {\footnotesize 2, \phantom{0}4} \\
		                                                 &                          &                          &     {\footnotesize 0.15} &                   {\footnotesize 10} \\ \cline{2-5}
		                                                 & \multirow{2}{*}{13\,828} & \multirow{3}{*}{0.00106} &     {\footnotesize 0.04} &                {\footnotesize 2, 10} \\
		                                                 &                          &                          &     {\footnotesize 0.08} &                {\footnotesize 2, 10} \\ \hline
		\multirow{3}{1.8cm}{\textbf{B2}}                 &                   5\,808 &                   0.0156 &                    0.005 &                                1, 10 \\ \cline{2-5}
		                                                 &                   4\,134 &                   0.0186 &                     0.01 &                                2, 10 \\ \cline{2-5}
		                                                 &                   8\,790 &                   0.0126 &                    0.013 &                                3, 10 \\ \hline
		\multirow{3}{1.8cm}{\textbf{B2S} \\ (scaled Bs)}  &                   5\,808 &                  0.00156 &                     0.05 &                                1, 10 \\ \cline{2-5}
		                                                 &                   4\,134 &                  0.00186 &                      0.1 &                                2, 10 \\ \cline{2-5}
		                                                 &                   8\,790 &                  0.00126 &                     0.13 &                                3, 10 \\ \hline
	\end{tabular}
	\caption{Different combinations of mesh sizes and flow velocities of the back-step simulations to create data. The number of obtained data points is given in column 2, which are sampled from the nonlinear iteration steps given in the last column. The dimensions of the four back-steps can be found in~\Cref{tab:backstepgeom}.
	\label{tab:data}
	}
\end{table}

We train our model based on simulation data from several boundary value problems defined on back-step geometries; see~\Cref{fig:CFLCOMSOL} for an exemplary back-step geometry. \AZ{The boundary conditions for the back-step flow are defined as follows: \AH{We impose a parabolic inflow velocity profile at the inlet $\Gamma_{in}$ and a constant pressure at the outlet $\Gamma_{out}$ on the right side of the geometry. For the other parts of the boundary, $\Gamma_{wall}$, a no-slip boundary condition is prescribed.}} 

As discussed in~\cref{sec:DataSets}, due to the locality of the input features, we expect that it suffices to use training data for a single type of geometry, as long as there is enough variation in the data with respect to the individual patches. In particular, we will consider two different base cases of back-step geometries, which we denote as B1 and B2. For both cases, we consider meshes with different levels of refinement and inflow velocities. This also results in flow fields with different Reynolds numbers. Moreover, we use data from different Newton iterations since, as can be seen in~\Cref{fig:b1}, the iterate change{\AH s} quite drastically between different Newton iterations. Additionally, we consider scaled variants of B1 and B2, which we denote as B1S and B2S, respectively. These are obtained by scaling the geometric dimensions of B1 and B2 by a factor of $0.1$ and \AH{increasing} the inflow velocity by a factor of $10$;  cf.~\Cref{tab:backstepgeom} for the dimensions of the four different geometries. As a result the Reynolds number remains the same, and we would expect the nonlinear convergence to be similar. By adding the scaled versions to the \AH{training configurations}, we try to force the model to learn the dependence of the nonlinearity on the Reynolds number. We summarize all configurations considered for the generation of training data in~\Cref{tab:data}.

In total, we obtain 79\,256 data points \AH{from all configurations listed in~\Cref{tab:data}. To arrive at a balanced distribution of the training data set, we perform the sampling such that we obtain roughly the same number of data points for each element size in the data set; that is, we select exactly 3\,500 data points from all configurations with the same maximum element size.} In total, we can categorize 14 different simulations based on the maximum element size, which gives a data set of \AZ{4}9\,000 data points. To make it easier to applied batch learning for the NN, 1\,000 data points are discarded at random. So, in total there are 48\,000 data points, of which 70\% \AH{are} used for training our network model, 15\% \AH{are} used for testing and 15\% \AH{are} used for validati\AH{on}. In~\cref{sec:results}, we will first discuss how this model performs on configurations \AH{with} back-step geometr\AH{y} in~\cref{sec:results:backstep}. Then, in~\cref{sec:results:couette,sec:results:obstacle}, we will also discuss the generalization to other configurations, that is, to Couette flow and flow around an obstacle, respectively.

%% file: Mainmatter/NetworkForCFLpred/backstep.tex
\begin{tikzpicture}
\draw[black,very thick,very thick] (0,0.7) -- (2,0.7) -- node[midway,left]{$\Gamma_{wall}$} (2,0) -- (8,0) -- (8,1.2) node[midway,right]{$\Gamma_{out}$} -- node[midway,above]{$\Gamma_{wall}$} (0,1.2) -- (0,0.7) node[midway,left]{$\Gamma_{in}$};



\draw[red,very thick,<->] (0,0.7) -- (0,1.2);
\draw[blue,very thick,<->] (2,0) -- (8,0);

\end{tikzpicture}

%% file: Mainmatter/NetworkForCFLpred/NetworkStructure.tex
\subsection{Network architecture and training}\label{sec:netstruc}

\begin{figure}[t]
	\centering
	\includegraphics[width = 0.45\textwidth]{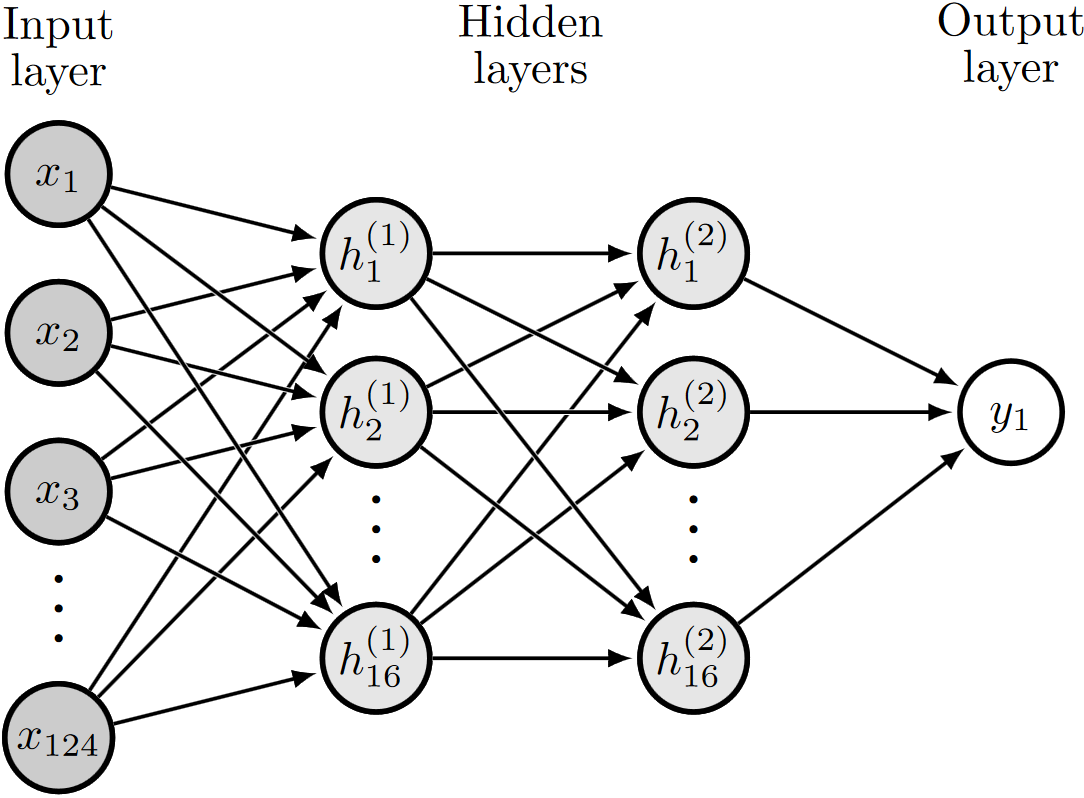}
	\caption{\AH{Exemplary NN} with 124 inputs, 2 hidden layers with 16 neurons and an output layer. The arrows are associated with weights in the weight matrices $W_i$\AH{; cf.~\cref{eq:nn}.}}
	\label{fig:nn1}
\end{figure}

The architecture of the network model employed in this work is depicted in~\Cref{fig:nn1}: the model consists of an input layer with 124 neurons, two hidden layers with \AH{16 neurons} each, and an output layer with the neuron representing prediction of the local optimal pseudo-time step; the composition of the input data has been discussed in~\cref{sec:DataSets}.

\begin{table}[t]
	\centering
	\begin{tabular}{|l|l|}
		\hline
		\textbf{hyper parameter} & \textbf{search space} \\
		\hline \hline
		\# hidden layers & $\left\lbrace 2,3,4,5 \right\rbrace$ \\
		\# neurons per hidden layer & $\left\lbrace 2^4,2^5,2^6,2^7,2^8 \right\rbrace$ \\
		\hline
	\end{tabular}
	\caption{Hyper parameters and search space for the grid search.
	\label{tab:grid_search}}
\end{table}

\begin{table}[t] \centering
	\begin{tabular}{|p{1.9cm}|p{3.3cm}|p{4.55cm}|}
		\hline
		\textbf{geometry type} & \textbf{inflow/wall \mbox{velocity} (\nicefrac{m}{s})}  & \textbf{\mbox{max.\ element size (m)}}             \\ \hline\hline
		\textbf{B1}            & 0.001,\,0.004,\,0.007, 0.01,\,0.012,\,0.015 & 0.0106\,:\,0.005\,:\,0.0256          \\ \hline
		\textbf{B1S}            & 0.01,\,0.04,\,0.07, 0.1,\,0.12,\,0.15       & 0.00106\,:\,0.0005\,:\,0.00256       \\ \hline
		\textbf{B2}            & 0.001\,:\,0.003\,:\,0.01                    & 0.0126,\,0.0156,\,0.0186, 0.0206     \\ \hline
		\textbf{B2S}            & 0.01\,:\,0.03\,:\,0.1                       & 0.00126,\,0.00156,\,0.00186, 0.00206 \\ \hline
		\textbf{C}            & 0.01,\,0.03,\,0.04, 0.05,\,0.07,\,0.1       & 0.014\,:\,0.002\,:\,0.022            \\ \hline
		\textbf{CS}            & 0.01,\,0.03,\,0.05                          & 0.0028\,:\,0.0004\,:\,0.0044         \\ \hline
	\end{tabular}
	\caption{Combinations of velocities and mesh sizes for simulations used to test the NN performance, with the back-step \AH{configurations B1--B2S} in~\Cref{tab:backstepgeom} and the Couette flow \AH{configurations} C1 \AH{and} C2 in~\Cref{tab:couettegeom}.
		\label{tab:valsim}
	}
\end{table}

We determined the network architecture as follows: We first performed a grid search with $6$-fold cross validation 
on the search space given in~\Cref{tab:grid_search} to determine three different network architectures with an acceptable average validation loss. \AH{N}ote that the loss is defined based on the prediction of the local optimal pseudo-time step\AH{; t}hat means a low validation loss does not necessarily result in a low solution residual when the network is used. 
\AH{From the grid search,} we determined one network with small, one with medium, and one with high networks capacity, that is, networks with:
\begin{itemize}
	\item 2 hidden layers and 16 neurons per layer,
	\item 3 hidden layers and 64 neurons per layer,
	\item 4 hidden layers and 256 neurons per layer.
\end{itemize}
Then, in order to find the best model for accelerating the Newton convergence, we compare them when applied to several back-step and Couette simulations, further explained in ~\cref{sec:results:backstep} and~\cref{sec:results:couette}. The combinations of the inflow or wall velocities with the maximum element size used for these simulations are given in~\Cref{tab:valsim}.

For each network, the required number of nonlinear iterations per simulation was compared. Surprisingly, the largest network now performed much worse, while the smallest network performed best, which could be explained by an overfitting behavior.

\begin{figure}[t]
	\centering
	\includegraphics[width=0.8\textwidth]{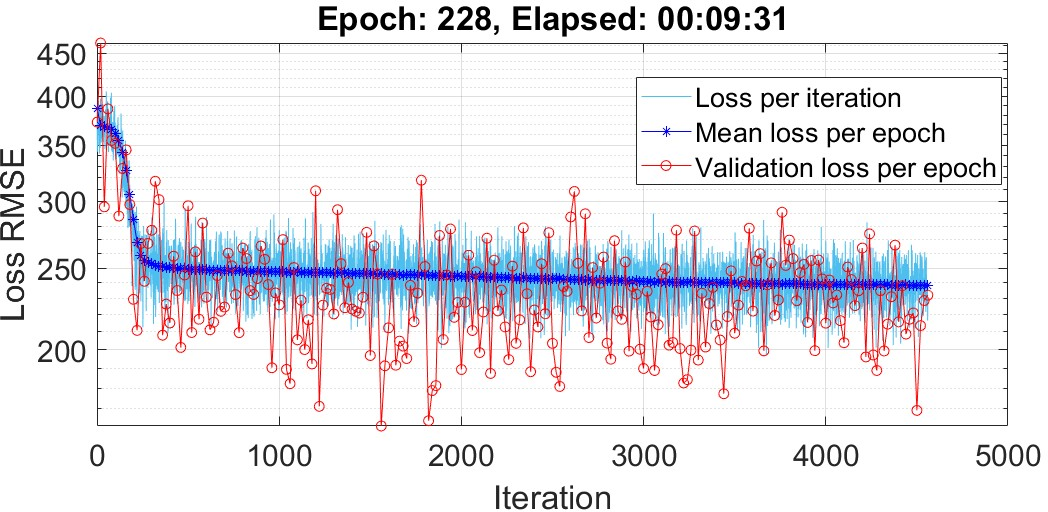}
	\caption{Training plot of network trained on the optimized pseudo-time step targets.}
	\label{fig:losscurve}
\end{figure}

Each network is \AH{trained} based on \AH{the} different back-step simulations given in~\Cref{tab:backstepgeom}. \AH{As mentioned before,} the local optimal pseudo-time steps obtained using the SNOPT optimizer \AH{are used} as \AH{the} reference. In particular, the network output is compared to these optimized local pseudo-time steps via the RMSE loss function~\cref{eq:rmseloss}. As the optimizer, we use stochastic gradient descent with adaptive moments (Adam)~\cite{kingma_adam_2017} with an initial learning rate of $0.001$. As the stopping criterion and regularization, we employ early stopping with a patience of $150$ epochs.

%% file: Mainmatter/Results/Results.tex
\section{Numerical results}
\label{sec:results}

\begin{figure}[t]
	\centering
	\includegraphics[width=0.55\textwidth]{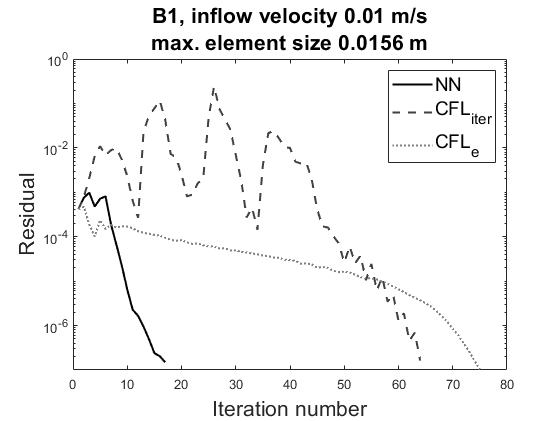}
	\caption{Convergence plot of a B1 simulation \AH{for all three strategies} for the local pseudo-time step \AH{under consideration: using the NN model as well as the strategies defined in~\cref{eq:CFLloc,eq:CFL2}.}
		\label{fig:convB1}}
\end{figure}

\begin{figure}[t]
	\centering
	
	\hfill
	\begin{minipage}{0.03\textwidth}
		\rotatebox{90}{\qquad \quad 2 \qquad it.}
	\end{minipage}
	\begin{minipage}{0.95\textwidth}
		\includegraphics[width=\textwidth]{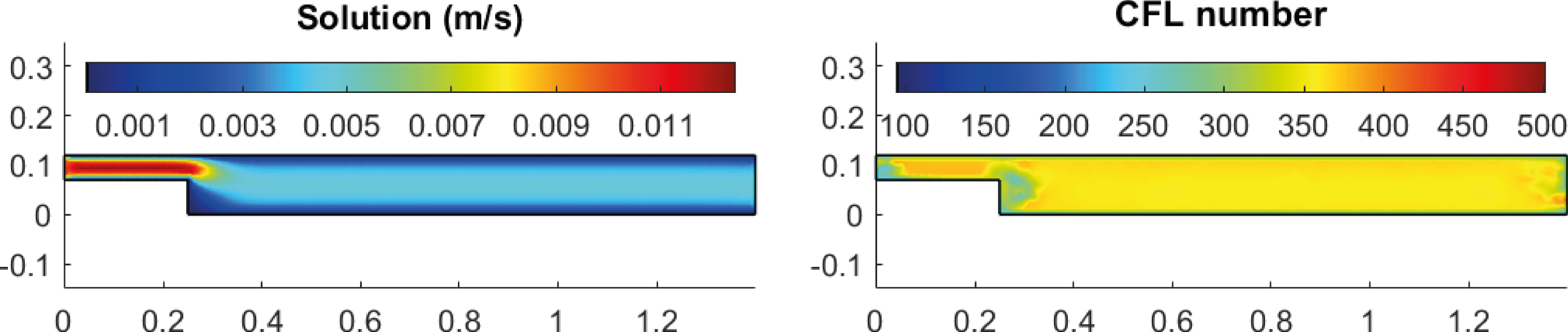}
	\end{minipage} \\[3.5ex]
	\hfill
	\begin{minipage}{0.03\textwidth}
		\rotatebox{90}{5}
	\end{minipage}
	\begin{minipage}{0.95\textwidth}
		\hfill
		\includegraphics[width=0.442\textwidth]{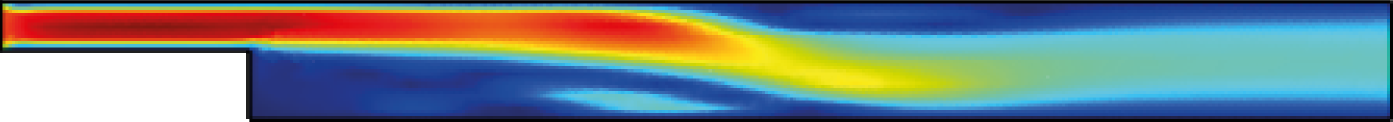}
		\hspace*{0.058\textwidth}
		\includegraphics[width=0.442\textwidth]{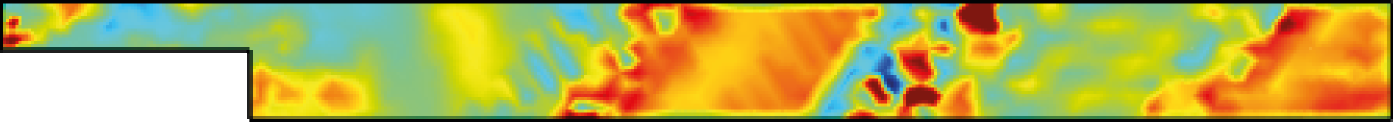}
	\end{minipage} \\[3.5ex]
	\hfill
	\begin{minipage}{0.03\textwidth}
		\rotatebox{90}{9}
	\end{minipage}
	\begin{minipage}{0.95\textwidth}
		\hfill
		\includegraphics[width=0.442\textwidth]{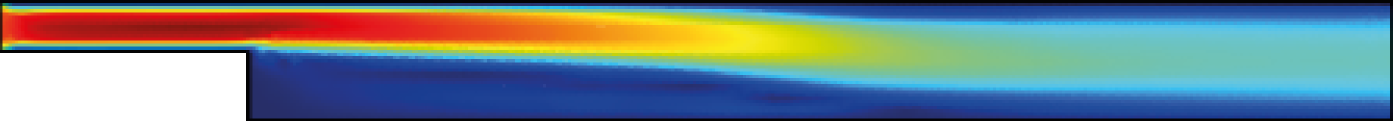}
		\hspace*{0.058\textwidth}
		\includegraphics[width=0.442\textwidth]{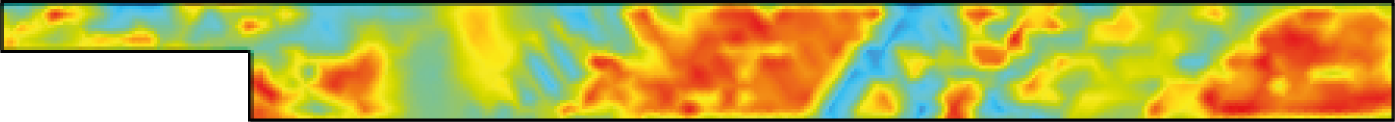}
	\end{minipage} \\[3.5ex]
	\hfill
	\begin{minipage}{0.03\textwidth}
		\rotatebox{90}{18}
	\end{minipage}
	\begin{minipage}{0.95\textwidth}
		\hfill
		\includegraphics[width=0.442\textwidth]{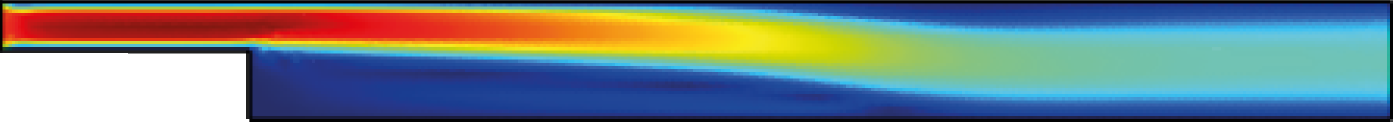}
		\hspace*{0.058\textwidth}
		\includegraphics[width=0.442\textwidth]{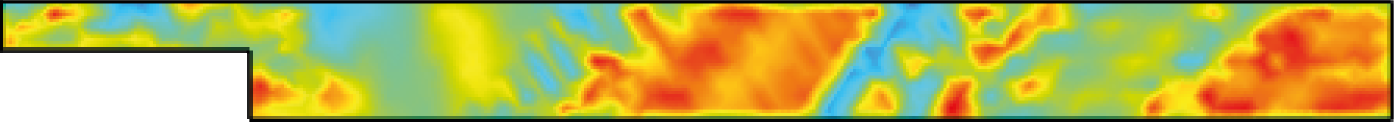}
	\end{minipage}
	\caption{Iterations of B1 with inflow velocity $0.001$\,\nicefrac{m}{s} and maximum element size $0.0156$\,m using $\text{NN}$, with velocities (left) and corresponding distribution of the predicted local pseudo-time step prediction (right). 
	}
	\label{fig:b1}
\end{figure}

\begin{figure}[t!]
	\begin{subfigure}{0.485\textwidth} \centering
		\input{Figures/Scatterplots/B1}
	\end{subfigure}
	\begin{subfigure}{0.485\textwidth} \centering
		\input{Figures/Scatterplots/B2}
	\end{subfigure}
	\begin{subfigure}{0.485\textwidth} \centering
		\input{Figures/Scatterplots/B3}
	\end{subfigure}
	\begin{subfigure}{0.485\textwidth} \centering
		\input{Figures/Scatterplots/B4}
	\end{subfigure}
	\begin{subfigure}{0.485\textwidth} \centering
		\input{Figures/Scatterplots/BM}
	\end{subfigure}
	\hspace{9pt}
	\begin{subfigure}{0.485\textwidth} \centering
		\input{Figures/Scatterplots/BR}
	\end{subfigure}
	\caption{The results of the back-steps for each method; the simulations are ordered from smallest required number of nonlinear iterations to largest \AH{separately for each approach}. \AH{The horizontal lines indicate the average nonlinear iteration counts for the different approaches.}
	}
	\label{fig:resultsB}
\end{figure}
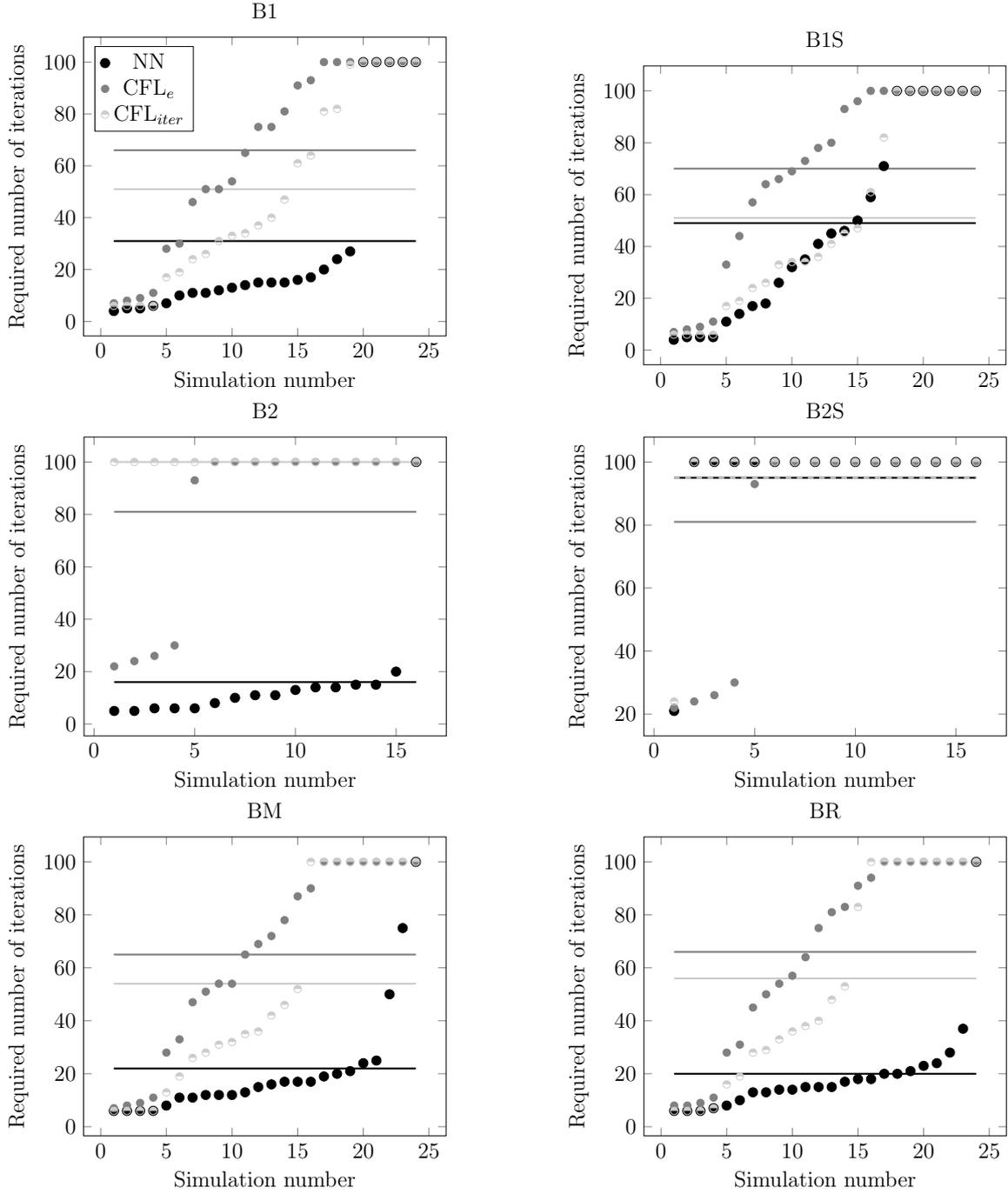

\begin{table}[t]
	\centering \small
	\begin{tabular}{|l|r|r|r|r|}
		\hline
		\textbf{Couette}  & \textbf{radius outer} &   \textbf{radius inner} &     \textbf{range wall} &               \textbf{range maximum} \\
		\textbf{geometry} &   \textbf{circle (m)} & \textbf{circle (m)} & \textbf{velocity (\nicefrac{m}{s})} &            \textbf{element size (m)} \\ \hline\hline
		\textbf{C}       &        0.4\phantom{0} &          0.2\phantom{0} &  0.01 -- 0.1\phantom{0} & 0.014\phantom{0} -- 0.022\phantom{0} \\ \hline
		\textbf{CS}       &                  0.08 &                    0.04 &            0.01 -- 0.05 &                     0.0028 -- 0.0044 \\ \hline
	\end{tabular}
	\caption{Geometrical dimensions of the two different Couette cylinders with corresponding ranges of velocities for the boundary conditions and maximum element sizes.}
	\label{tab:couettegeom}
\end{table}

\begin{figure}[t]
	\centering
	
	\hspace*{0.01\textwidth} 2 \hspace*{0.19\textwidth} 5 \hspace*{0.19\textwidth} 9 \hspace*{0.185\textwidth} 18 \hspace*{0.02\textwidth} \hfill \\
	\includegraphics[height=0.225\textwidth]{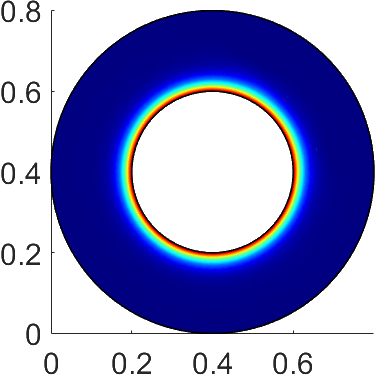}
	\hspace*{0.01\textwidth}
	\includegraphics[height=0.225\textwidth]{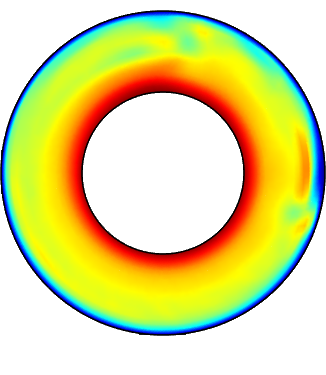}
	\hspace*{0.01\textwidth}
	\includegraphics[height=0.225\textwidth]{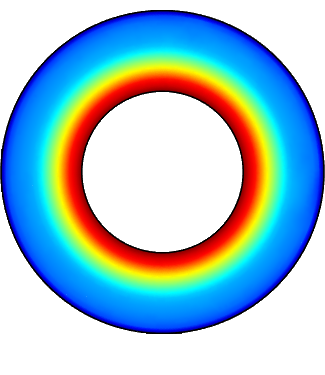}
	\hspace*{0.01\textwidth}
	\includegraphics[height=0.225\textwidth]{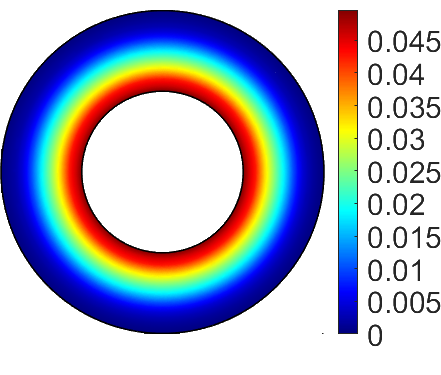} \\[1ex]
	\includegraphics[height=0.225\textwidth]{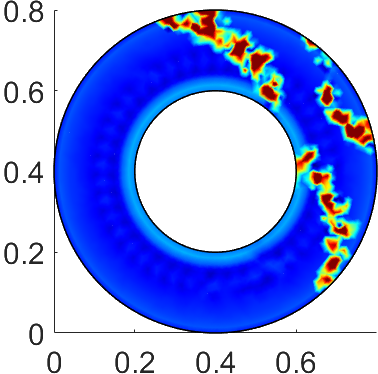}
	\hspace*{0.01\textwidth}
	\includegraphics[height=0.225\textwidth]{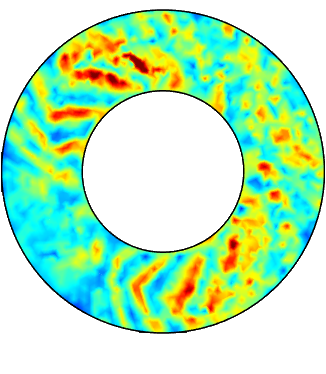}
	\hspace*{0.01\textwidth}
	\includegraphics[height=0.225\textwidth]{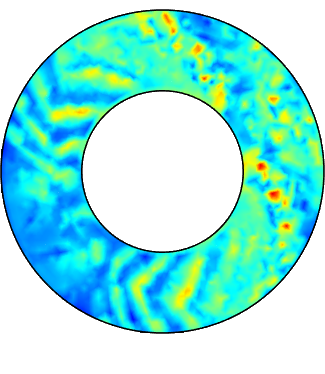}
	\hspace*{0.01\textwidth}
	\includegraphics[height=0.225\textwidth]{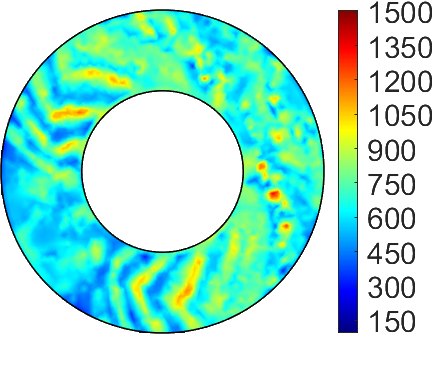}
	\caption{Iterations of C with wall velocity 0.05 \nicefrac{m}{s} and maximum element size of 0.016 m using $\text{NN}$, with velocity (above) and local pseudo-time step prediction (below).\label{fig:C}
	}
\end{figure}

\begin{figure}[t]
	\centering
	\includegraphics[width=0.55\textwidth]{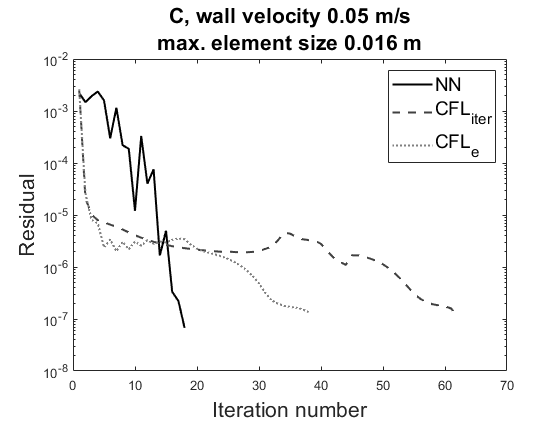}
	\caption{Convergence plot of a C simulation \AH{for all three strategies} for the local pseudo-time step \AH{under consideration: using the NN model as well as the strategies defined in~\cref{eq:CFLloc,eq:CFL2}.}
		\label{fig:convC}} 
\end{figure}

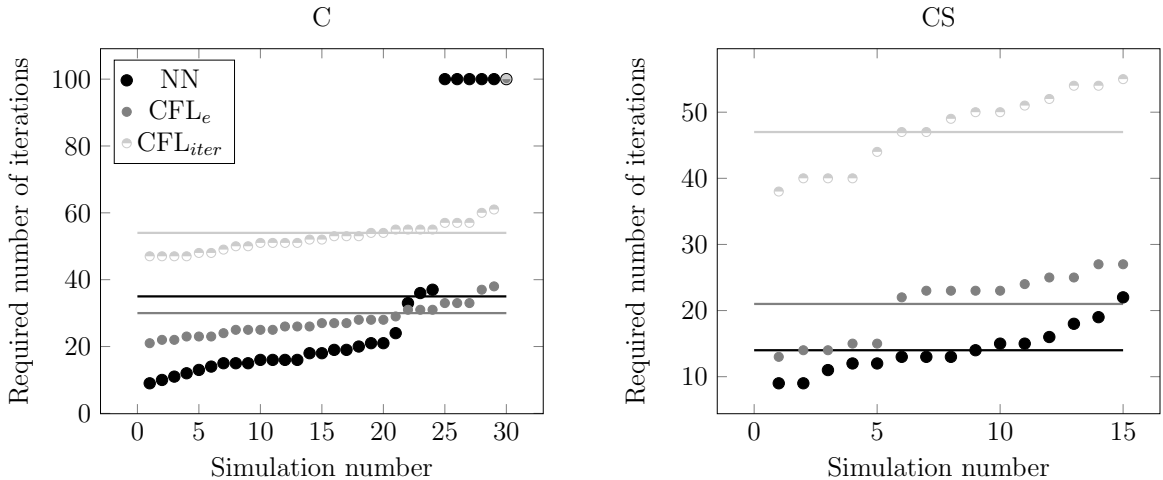
\begin{figure}[t]
	\begin{subfigure}{0.485\textwidth} \centering
		\input{Figures/Scatterplots/C1}
	\end{subfigure}
	\begin{subfigure}{0.485\textwidth} \centering
		\input{Figures/Scatterplots/C2}
	\end{subfigure}
	\caption{The results of the Couette flow for each method; the simulations are ordered from smallest required number of nonlinear iterations to largest \AH{separately for each approach. \AH{The horizontal lines indicate the average nonlinear iteration counts for the different approaches.}
			\label{fig:resultsC}}
	}
\end{figure}

In this section, we discuss numerical results for our NN-enhanced pseudo-time stepping approach, which we have \AH{introduced} in the previous~\cref{sec:NN}. As discussed in~\cref{sec:traindata}, we have trained the model only for back-step geometries. Our model is designed such that only local features are used as input in order to enable generalizability. First, we test our model on the \textit{back-step flow} \AH{training} cases, with \AH{varying} mesh sizes and inflow velocities within the range of the training data; cf.~\cref{sec:results:backstep}. \AH{In addition to that, we investigate the performance on flipped and rotated cases.} Then, in order to \AH{further} test the generalization properties of the model, we test our model on \textit{Couette cylinder flows} in~\cref{sec:results:couette} as well as for \textit{flow around an obstacle} cases in~\cref{sec:results:obstacle}. For the Couette cases, we vary the geometry, the mesh size, and the velocity of the inner rotating wall, and for the flow around an obstacle, back-step and Couette cylinder geometries with different kinds of obstacles are considered.

In particular, we compare the convergence of our hybrid approach with the convergence of two classical strategies for the choice of the local pseudo-time step described in~\cref{sec:local_cfl}, that is, the strategy denoted by CFL$_{iter}$, as given by~\cref{eq:deltaT,eq:CFLloc}, and the strategy denoted by CFL$_e$, as given by~\cref{eq:deltaT,eq:CFL2}. Furthermore, \AH{in the appendix, the convergence of the network is also compared with two variants of Newton's method, that is, the standard Newton method with a constant damping factor and a variant with an adaptive choice of the damping factor; see~\cref{sec:appendix} for more details}. The $L^2$ norm of the "reconstructed" residual $\tilde{R}$
\begin{equation} \label{eq:E}
	\|  \tilde{R} \|_{L^2(\Omega)} 
	=
	\sqrt{
		\int_\Omega
		\tilde{R}_u^2 +
		\tilde{R}_v^2 +
		\tilde{R}_p^2
		\, dx
	}
\end{equation}
is used as the stopping criterion. Here, $(\tilde{R}_u,\tilde{R}_v)$ and $\tilde{R}_p$ are the residual velocity and pressure fields in $X_h$ reconstructed by using the components of the residual vector $R=(R_u,R_v,R_p)$, as defined in~\cref{eq:R}, as coefficients w.r.t.\ the chosen basis of $X_h$.

Note that COMSOL uses a scaled version of the residual~\cref{eq:E} in the stopping criterion, such that the tolerance may differ depending on the problem configuration. However, when reporting the convergence results, we make sure that the same tolerance for the unscaled residual~\cref{eq:E} is used to determine convergence for the different approaches.

\AH{All simulations have been performed using COMSOL Multiphysics\textsuperscript{\tiny\textregistered} version 6.1.}

\input{Mainmatter/Results/BackStep}

\input{Mainmatter/Results/Couette}

\input{Mainmatter/Results/Obstacle}

%% file: Figures/Scatterplots/B1.tex
\begin{tikzpicture}[scale=0.8]
\begin{axis}[
    title = B1,
    legend pos = north west,
    xlabel={Simulation number},
    ylabel={Required number of iterations}
    ]
    \addplot [
        scatter,
        only marks,
        point meta=explicit symbolic,
        scatter/classes={
            a={mark size=2.5pt,mark=*,black},
            b={mark size=2pt, mark=*,black!50},
            c={mark size=2pt, mark=myhalfcircle2,black!20}
        },
    ] table [meta=label] {
        x      y        label
        1 4 a
2 5 a
3 5 a
4 6 a
5 7 a
6 10 a
7 11 a
8 11 a
9 12 a
10 13 a
11 14 a
12 15 a
13 15 a
14 15 a
15 16 a
16 17 a
17 20 a
18 24 a
19 27 a
20 100 a
21 100 a
22 100 a
23 100 a
24 100 a
1 7 b
2 8 b
3 9 b
4 11 b
5 28 b
6 30 b
7 46 b
8 51 b
9 51 b
10 54 b
11 65 b
12 75 b
13 75 b
14 81 b
15 91 b
16 93 b
17 100 b
18 100 b
19 100 b
20 100 b
21 100 b
22 100 b
23 100 b
24 100 b
1 6 c
2 6 c
3 6 c
4 6 c
5 17 c
6 19 c
7 24 c
8 26 c
9 31 c
10 33 c
11 34 c
12 37 c
13 40 c
14 47 c
15 61 c
16 64 c
17 81 c
18 82 c
19 99 c
20 100 c
21 100 c
22 100 c
23 100 c
24 100 c
    };    
    \addplot[line width=1pt, color = black, 
] coordinates {(1,31) (24,31)};
\addplot[line width=1pt, color = black!50, 
] coordinates {(1,66) (24,66)};
\addplot[line width=1pt, color = black!20, 
] coordinates {(1,51) (24,51)};
\legend{NN,CFL$_e$,CFL$_{iter}$}
\end{axis}
\end{tikzpicture}

%% file: Figures/Scatterplots/B2.tex
\begin{tikzpicture}[scale=0.8]
\begin{axis}[
    title = B1S,
    legend pos = north west,
    ylabel={Required number of iterations}
    ]
    \addplot [
        scatter,
        only marks,
        point meta=explicit symbolic,
        scatter/classes={
            a={mark size=2.5pt,mark=*,black},
            b={mark size=2pt, mark=*,black!50},
            c={mark size=2pt, mark=myhalfcircle2,black!20}
        },
    ] table [meta=label] {
        x y label
        1 4 a
2 5 a
3 5 a
4 5 a
5 11 a
6 14 a
7 17 a
8 18 a
9 26 a
10 32 a
11 35 a
12 41 a
13 45 a
14 46 a
15 50 a
16 59 a
17 71 a
18 100 a
19 100 a
20 100 a
21 100 a
22 100 a
23 100 a
24 100 a
1 7 b
2 8 b
3 9 b
4 11 b
5 33 b
6 44 b
7 57 b
8 64 b
9 66 b
10 69 b
11 73 b
12 78 b
13 80 b
14 93 b
15 96 b
16 100 b
17 100 b
18 100 b
19 100 b
20 100 b
21 100 b
22 100 b
23 100 b
24 100 b
1 6 c
2 6 c
3 6 c
4 6 c
5 17 c
6 19 c
7 24 c
8 26 c
9 33 c
10 34 c
11 34 c
12 36 c
13 41 c
14 45 c
15 47 c
16 61 c
17 82 c
18 100 c
19 100 c
20 100 c
21 100 c
22 100 c
23 100 c
24 100 c
    };   
\addplot[line width=1pt, color = black, 
] coordinates {(1,49) (24,49)};
\addplot[line width=1pt, color = black!50, 
] coordinates {(1,70) (24,70)};
\addplot[line width=1pt, color = black!20, 
] coordinates {(1,51) (24,51)};
\end{axis}
\end{tikzpicture}

%% file: Figures/Scatterplots/B3.tex
\begin{tikzpicture}[scale=0.8]
\begin{axis}[
    title = B2,
    legend pos = north west,
    xlabel={Simulation number},
    ylabel={Required number of iterations}
    ]
    \addplot [
        scatter,
        only marks,
        point meta=explicit symbolic,
        scatter/classes={
            a={mark size=2.5pt,mark=*,black},
            b={mark size=2pt, mark=*,black!50},
            c={mark size=2pt, mark=myhalfcircle2,black!20}
        },
    ] table [meta=label] {
        x y label
        1 5 a
2 5 a
3 6 a
4 6 a
5 6 a
6 8 a
7 10 a
8 11 a
9 11 a
10 13 a
11 14 a
12 14 a
13 15 a
14 15 a
15 20 a
16 100 a
1 22 b
2 24 b
3 26 b
4 30 b
5 93 b
6 100 b
7 100 b
8 100 b
9 100 b
10 100 b
11 100 b
12 100 b
13 100 b
14 100 b
15 100 b
16 100 b
1 100 c
2 100 c
3 100 c
4 100 c
5 100 c
6 100 c
7 100 c
8 100 c
9 100 c
10 100 c
11 100 c
12 100 c
13 100 c
14 100 c
15 100 c
16 100 c
    };    
    \addplot[line width=1pt, color = black, 
] coordinates {(1,16) (16,16)};
\addplot[line width=1pt, color = black!50, 
] coordinates {(1,81) (16,81)};
\addplot[line width=1pt, color = black!20, 
] coordinates {(1,100) (16,100)};
\end{axis}
\end{tikzpicture}

%% file: Figures/Scatterplots/B4.tex
\begin{tikzpicture}[scale=0.8]
\begin{axis}[
    title = B2S,
    legend pos = south east,
    xlabel={Simulation number},
    ylabel={Required number of iterations}
    ]
    \addplot [
        scatter,
        only marks,
        point meta=explicit symbolic,
        scatter/classes={
            a={mark size=2.5pt,mark=*,black},
            b={mark size=2pt, mark=*,black!50},
            c={mark size=2pt, mark=myhalfcircle2,black!20}
        },
    ] table [meta=label] {
        x y label
1 21 a
2 100 a
3 100 a
4 100 a
5 100 a
6 100 a
7 100 a
8 100 a
9 100 a
10 100 a
11 100 a
12 100 a
13 100 a
14 100 a
15 100 a
16 100 a
1 22 b
2 24 b
3 26 b
4 30 b
5 93 b
6 100 b
7 100 b
8 100 b
9 100 b
10 100 b
11 100 b
12 100 b
13 100 b
14 100 b
15 100 b
16 100 b
1 24 c
2 100 c
3 100 c
4 100 c
5 100 c
6 100 c
7 100 c
8 100 c
9 100 c
10 100 c
11 100 c
12 100 c
13 100 c
14 100 c
15 100 c
16 100 c
    };   
    \addplot[line width=1pt, color = black, 
] coordinates {(1,95) (16,95)};
\addplot[line width=1pt, color = black!50, 
] coordinates {(1,81) (16,81)};
\addplot[line width=1pt, color = black!20, 
style = dashed
] coordinates {(1,95) (16,95)};
\end{axis}
\end{tikzpicture}

%% file: Figures/Scatterplots/BM.tex
\begin{tikzpicture}[scale=0.8]
\begin{axis}[
    title = BM,
    legend pos = north west,
    xlabel={Simulation number},
    ylabel={Required number of iterations}
    ]
    \addplot [
        scatter,
        only marks,
        point meta=explicit symbolic,
        scatter/classes={
            a={mark size=2.5pt,mark=*,black},
            b={mark size=2pt, mark=*,black!50},
            c={mark size=2pt, mark=myhalfcircle2,black!20}
        },
    ] table [meta=label] {
        x y label
        1 6 a
2 6 a
3 6 a
4 6 a
5 8 a
6 11 a
7 11 a
8 12 a
9 12 a
10 12 a
11 13 a
12 15 a
13 16 a
14 17 a
15 17 a
16 17 a
17 19 a
18 20 a
19 21 a
20 24 a
21 25 a
22 50 a
23 75 a
24 100 a
1 7 b
2 8 b
3 9 b
4 11 b
5 28 b
6 33 b
7 47 b
8 51 b
9 54 b
10 54 b
11 65 b
12 69 b
13 72 b
14 78 b
15 87 b
16 90 b
17 100 b
18 100 b
19 100 b
20 100 b
21 100 b
22 100 b
23 100 b
24 100 b
1 6 c
2 6 c
3 6 c
4 6 c
5 13 c
6 19 c
7 26 c
8 28 c
9 31 c
10 32 c
11 35 c
12 36 c
13 42 c
14 46 c
15 52 c
16 100 c
17 100 c
18 100 c
19 100 c
20 100 c
21 100 c
22 100 c
23 100 c
24 100 c
 };    
 \addplot[line width=1pt, color = black, 
] coordinates {(1,22) (24,22)};
\addplot[line width=1pt, color = black!50, 
] coordinates {(1,65) (24,65)};
\addplot[line width=1pt, color = black!20, 
] coordinates {(1,54) (24,54)};
\end{axis}
\end{tikzpicture}

%% file: Figures/Scatterplots/BR.tex
\begin{tikzpicture}[scale=0.8]
\begin{axis}[
    title = BR,
    legend pos = north west,
    xlabel={Simulation number},
    ylabel={Required number of iterations}
    ]
    \addplot [
        scatter,
        only marks,
        point meta=explicit symbolic,
        scatter/classes={
            a={mark size=2.5pt,mark=*,black},
            b={mark size=2pt, mark=*,black!50},
            c={mark size=2pt, mark=myhalfcircle2,black!20}
        },
    ] table [meta=label] {
        x y label
        1 6 a
2 6 a
3 6 a
4 7 a
5 8 a
6 10 a
7 13 a
8 13 a
9 14 a
10 14 a
11 15 a
12 15 a
13 15 a
14 17 a
15 18 a
16 18 a
17 20 a
18 20 a
19 21 a
20 23 a
21 24 a
22 28 a
23 37 a
24 100 a
1 8 b
2 8 b
3 9 b
4 11 b
5 28 b
6 31 b
7 45 b
8 50 b
9 54 b
10 57 b
11 64 b
12 75 b
13 81 b
14 83 b
15 91 b
16 94 b
17 100 b
18 100 b
19 100 b
20 100 b
21 100 b
22 100 b
23 100 b
24 100 b
1 6 c
2 6 c
3 6 c
4 7 c
5 16 c
6 19 c
7 28 c
8 29 c
9 33 c
10 36 c
11 38 c
12 40 c
13 48 c
14 53 c
15 83 c
16 100 c
17 100 c
18 100 c
19 100 c
20 100 c
21 100 c
22 100 c
23 100 c
24 100 c
 };    
 \addplot[line width=1pt, color = black, 
] coordinates {(1,20) (24,20)};
\addplot[line width=1pt, color = black!50, 
] coordinates {(1,66) (24,66)};
\addplot[line width=1pt, color = black!20, 
] coordinates {(1,56) (24,56)};
\end{axis}
\end{tikzpicture}

%% file: Figures/Scatterplots/C1.tex
\begin{tikzpicture}[scale=0.85]
\begin{axis}[
    title = C,
    legend pos = north west,
    xlabel={Simulation number},
    ylabel={Required number of iterations}
    ]
    \addplot [
        scatter,
        only marks,
        point meta=explicit symbolic,
        scatter/classes={
            a={mark size=2.5pt,mark=*,black},
            b={mark size=2pt, mark=*,black!50},
            c={mark size=2pt, mark=myhalfcircle2,black!20}
        },
    ] table [meta=label] {
        x y label
1 9 a
2 10 a
3 11 a
4 12 a
5 13 a
6 14 a
7 15 a
8 15 a
9 15 a
10 16 a
11 16 a
12 16 a
13 16 a
14 18 a
15 18 a
16 19 a
17 19 a
18 20 a
19 21 a
20 21 a
21 24 a
22 33 a
23 36 a
24 37 a
25 100 a
26 100 a
27 100 a
28 100 a
29 100 a
30 100 a
1 21 b
2 22 b
3 22 b
4 23 b
5 23 b
6 23 b
7 24 b
8 25 b
9 25 b
10 25 b
11 25 b
12 26 b
13 26 b
14 26 b
15 27 b
16 27 b
17 27 b
18 28 b
19 28 b
20 28 b
21 29 b
22 31 b
23 31 b
24 31 b
25 33 b
26 33 b
27 33 b
28 37 b
29 38 b
30 100 b
1 47 c
2 47 c
3 47 c
4 47 c
5 48 c
6 48 c
7 49 c
8 50 c
9 50 c
10 51 c
11 51 c
12 51 c
13 51 c
14 52 c
15 52 c
16 53 c
17 53 c
18 53 c
19 54 c
20 54 c
21 55 c
22 55 c
23 55 c
24 55 c
25 57 c
26 57 c
27 57 c
28 60 c
29 61 c
30 100 c
 };   
 \addplot[line width=1pt, color = black, 
] coordinates {(0,35) (30,35)};
\addplot[line width=1pt, color = black!50, 
] coordinates {(0,30) (30,30)};
\addplot[line width=1pt, color = black!20, 
] coordinates {(0,54) (30,54)};
\legend{NN,CFL$_e$,CFL$_{iter}$}
\end{axis}
\end{tikzpicture}

%% file: Figures/Scatterplots/C2.tex
\begin{tikzpicture}[scale=0.85]
\begin{axis}[
    title = CS,
    legend pos = north west,
    xlabel={Simulation number},
    ylabel={Required number of iterations}
    ]
    \addplot [
        scatter,
        only marks,
        point meta=explicit symbolic,
        scatter/classes={
            a={mark size=2.5pt,mark=*,black},
            b={mark size=2pt, mark=*,black!50},
            c={mark size=2pt, mark=myhalfcircle2,black!20}
        },
    ] table [meta=label] {
       x y label
1 9 a
2 9 a
3 11 a
4 12 a
5 12 a
6 13 a
7 13 a
8 13 a
9 14 a
10 15 a
11 15 a
12 16 a
13 18 a
14 19 a
15 22 a
1 13 b
2 14 b
3 14 b
4 15 b
5 15 b
6 22 b
7 23 b
8 23 b
9 23 b
10 23 b
11 24 b
12 25 b
13 25 b
14 27 b
15 27 b
1 38 c
2 40 c
3 40 c
4 40 c
5 44 c
6 47 c
7 47 c
8 49 c
9 50 c
10 50 c
11 51 c
12 52 c
13 54 c
14 54 c
15 55 c
 };    
 \addplot[line width=1pt, color = black, 
] coordinates {(0,14) (15,14)};
\addplot[line width=1pt, color = black!50, 
] coordinates {(0,21) (15,21)};
\addplot[line width=1pt, color = black!20, 
] coordinates {(0,47) (15,47)};
\end{axis}
\end{tikzpicture}

%% file: Mainmatter/Results/BackStep.tex
\subsection{Back-step geometries} \label{sec:results:backstep}

In this subsection, we discuss the performance of our model on the same four back-step geometries which we also used for the training; cf.~\cref{sec:traindata,tab:backstepgeom}. In particular, we investigate if the network can accelerate the convergence when varying inflow velocity and mesh size on those configurations; cf.~\Cref{tab:valsim}. As a result of those variations, a total of $80$ different simulations for the back-step geometries are analyzed.

Exemplary, we plot the convergence history for the back-step geometry B1 with inflow velocity $0.001$\,\nicefrac{m}{s} and a maximum element size of $0.0156$\,m in~\Cref{fig:convB1}; in particular, we compare the local optimal pseudo-time step predicted by the network model, which we denote as $\text{NN}$, against the $\text{CFL}_{iter}$ and CFL$_e$ approaches based on the iteration count and the evolution of the residual norm. It can be observed that the network speeds up the convergence significantly compared with the two reference approaches; moreover, whereas the residual norm oscillates strongly for $\text{CFL}_{iter}$, the convergence is almost monotonous for the network and CFL$_e$ approaches. \AH{The qualitative behavior observed in~\Cref{fig:convB1}, where the classical choices for the CFL number yield slow convergence or even strong oscillations in the residual, is not an exception; similar convergence behavior can also be observed for Couette flow and flow around an obstacle configurations; cf.~\cref{sec:results:couette,sec:results:obstacle}.} Furthermore, in~\Cref{fig:b1}, we depict four iterates during the Newton iteration using the network approach as well as the corresponding predicted local optimal pseudo-time step. It can be observed that both the iterate and predicted local optimal pseudo-time step still vary during the first iterations and both become stationary in the later iterations.

In\AH{~\Cref{fig:resultsB}}, we compare the performance with respect to the numbers of Newton iterations resulting from the different choices $\text{NN}$, $\text{CFL}_{iter}$, and CFL$_e$ for the pseudo-time step on different back-step configurations; cf.~\cref{sec:traindata,tab:backstepgeom}. Furthermore, as a first investigation of the generalizability, we additionally consider cases resulting from mirroring (BM) or anti-clockwise rotating by 90 degrees (BR) the B1 cases. 

For most cases from our training data B1--B2S, the NN model performs better than or equally well as the best choice of $\text{CFL}_{iter}$ and CFL$_e$, as well as both Newton methods\AH{; see~\cref{sec:appendix} for the results for the Newton methods}. 
\AH{This can also be seen in the average number of iterations for the configurations B1, B2, BM, and BR; see~\Cref{fig:resultsB,fig:resultsBSO}. For B1S, the network performs equally well as CFL$_{iter}$.}

In particular, for cases with lower velocities (B1 and B2), the network performs best. For cases with higher velocities (B1S and B2S) but the same Reynolds numbers, the network model performs clearly worse \AH{compared to the low velocity cases}. 
\AH{Moreover}, it \AH{could} be observed that, for B1S and B2S, \AH{convergence was generally less robust; in particular, we observed a relatively high ratio of cases for which none of the pseudo-time step approaches converged}. Since the nonlinearity should be mostly determined by the Reynolds number rather than the magnitude of the velocities, we would like to investigate in future research whether a scaling/normalization of the network inputs can improve the performance. \AZ{The same behavior can also be observed for \AH{the Newton methods,} although these methods converge in more cases\AH{; cf.~\cref{sec:appendix}.}}

Even though the mirrored and rotated configurations, of course, yield the same but mirrored respectively rotated results, we observe slight variations in the performance of the network model. This shows a slight overfitting with respect to the B1 configurations as shown exemplarily in~\Cref{fig:b1}, where, for instance, the velocities in positive $x$ direction are dominant; this clearly changes when the geometry is mirrored or rotated. However, the performance of the network model is only slightly worse compared to the B1 cases, which shows that the overfitting is not severe with respect to the prediction of suitable local pseudo-time steps. \AH{In the future, we may consider group invariant networks~\cite{cohen_group_2016} to prevent these geometric overfitting effects.}

%% file: Mainmatter/Results/Couette.tex
\subsection{Couette} \label{sec:results:couette}

Next, we discuss the generalization of our hybrid globalization approach to a different type of geometries, that is, Couette cylinders. We consider two different geometries and a range of configurations resulting from varying the wall rotation velocity and refinement of the mesh; cf.~\cref{tab:couettegeom}. In total, we obtain 45 different configurations of Couette simulations, and in~\cref{tab:valsim}, \AH{the corresponding} boundary conditions and maximum element sizes \AH{for the} back-step configurations previously used for training the NN model \AH{can be found}.

Again, we present the iterates and corresponding predicted local pseudo-time steps during the Newton iteration for one exemplary case with C geometry; see~\cref{fig:C}. As for the back-step configuration in~\cref{fig:b1}, we observe that both the solution and local pseudo-time step prediction approach a stationary distribution. For this configuration, we observe a significant speed-up of Newton convergence when using the local pseudo-time step prediction by our NN\AH{, which has only been trained on back-step geometries}; cf.~\cref{fig:convC}.

\AH{\Cref{fig:resultsC,fig:resultsCSO} confirm} that the use of our NN-based approach is beneficial for most of the Couette flow cases. First of all, we notice that, for the Couette flow simulations considered, all choices for the local pseudo-time step yield convergence. However, the NN-based approach accelerates convergence in 70\,\% of the C cases and even \AZ{93.3}\,\% of the CS cases. The fact that the network seems to perform better for CS than for C cases, compared with the $\text{CFL}_{iter}$ and CFL$_e$ choices, is remarkable because the CS cases contain higher velocities and smaller mesh element sizes compared to C. For the back-step geometries, we observed the opposite; cf.~\cref{fig:resultsB}, where the network performed worse on B1S and B2S than B1 and B2 configurations. 

\AH{Also when taking the two variants of Newton's method into account in the results in~\cref{sec:appendix}, the NN approach remains competitive. In particular, the NN is on average close to the best method (for the C configurations) or as good as the best method (for the CS configurations).}

%% file: Mainmatter/Results/Obstacle.tex
\begin{figure}[t]
    \begin{subfigure}[t]{0.29\textwidth}
        \centering 
        \includegraphics[width=\textwidth]{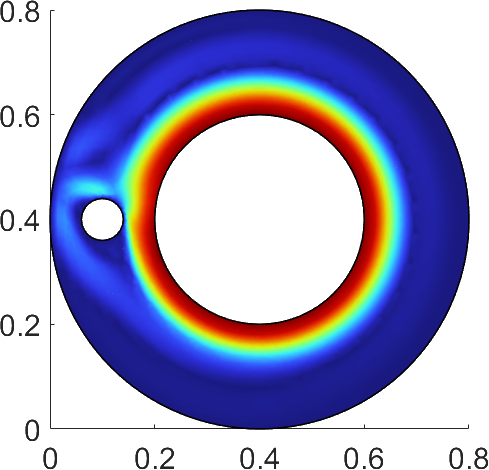}
        \caption{Circular obstacle.}
        \label{fig:couettecirc}
    \end{subfigure}
    \hfill 
    \begin{subfigure}[t]{0.29\textwidth}
        \centering 
        \includegraphics[width=\textwidth]{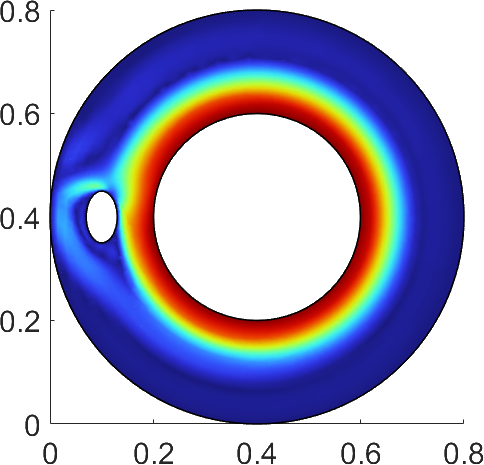}
        \caption{Elliptic obstacle with a-semiaxis $0.03$\,m and b-semiaxis $0.05$\,m.}
        \label{fig:couetteell}
    \end{subfigure}
    \hfill
    \begin{subfigure}[t]{0.36\textwidth}
        \centering 
        \includegraphics[width=\textwidth]{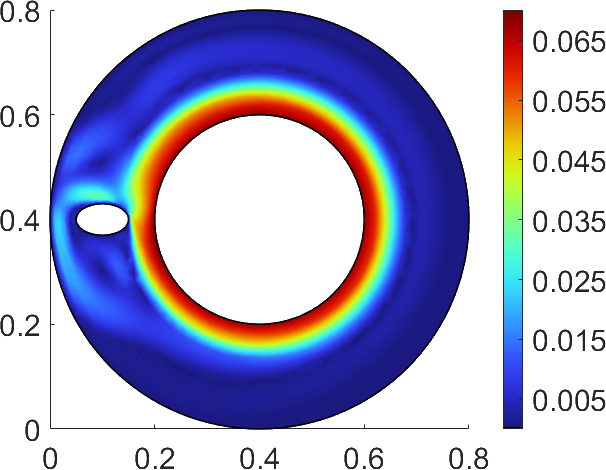}
        \caption{Elliptic obstacle with a-semiaxis $0.05$\,m and b-semiaxis $0.03$\,m.}
        \label{fig:couettecel2}
    \end{subfigure}

	\vspace*{-0.5cm}

    \caption{Velocity solutions for Couette flow with different kinds of obstacles. The wall rotation is 0.07 \nicefrac{m}{s} and the maximum element size is 0.02 m. The center of the obstacle is x = 0.1 m and y = 0.4 m.}\label{fig:couetteobst}
\end{figure}

\begin{figure}[t]
	\centering
	
	\hspace*{0.01\textwidth} 2 \hspace*{0.19\textwidth} 5 \hspace*{0.19\textwidth} 9 \hspace*{0.185\textwidth} 24 \hspace*{0.02\textwidth} \hfill \\
	\includegraphics[height=0.225\textwidth]{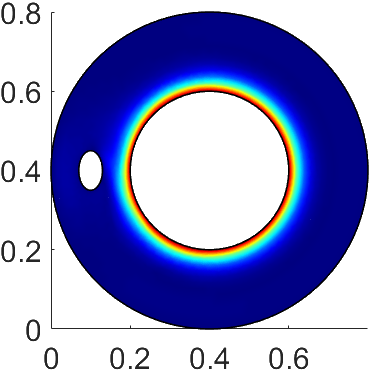}
	\hspace*{0.01\textwidth}
	\includegraphics[height=0.225\textwidth]{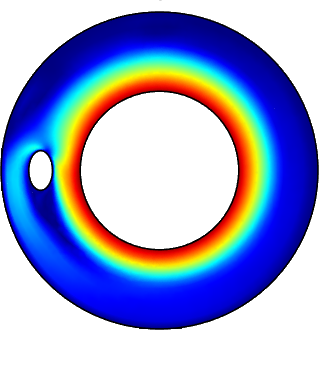}
	\hspace*{0.01\textwidth}
	\includegraphics[height=0.225\textwidth]{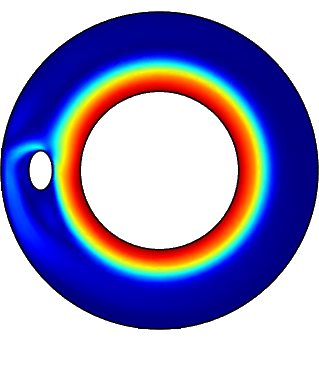}
	\hspace*{0.01\textwidth}
	\includegraphics[height=0.225\textwidth]{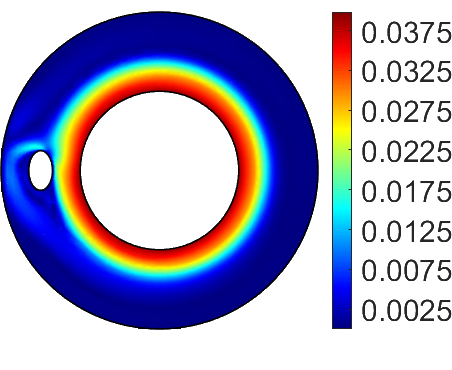} \\[1ex]
	\includegraphics[height=0.225\textwidth]{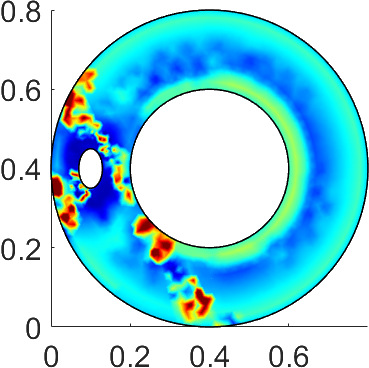}
	\hspace*{0.01\textwidth}
	\includegraphics[height=0.225\textwidth]{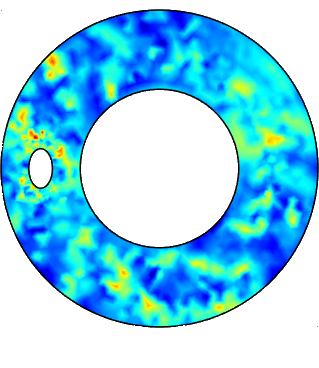}
	\hspace*{0.01\textwidth}
	\includegraphics[height=0.225\textwidth]{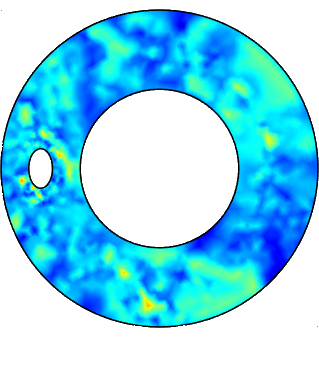}
	\hspace*{0.01\textwidth}
	\includegraphics[height=0.225\textwidth]{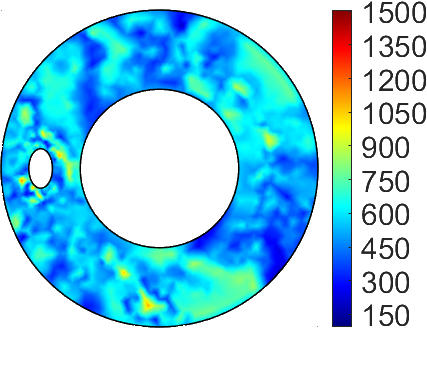}
	\caption{Iterations of CO with wall velocity 0.04 \nicefrac{m}{s} and maximum element size of 0.022 m using $\text{NN}$, with velocity (above) and local pseudo-time step prediction (below). The obstacle is an ellipse with a-semiaxis 0.03 m and b-semiaxis 0.05 m, with center x = 0.1 m and y = 0.4 m.
	}\label{fig:co}
\end{figure}

\begin{figure}[t]
	\centering
	\includegraphics[width=0.55\textwidth]{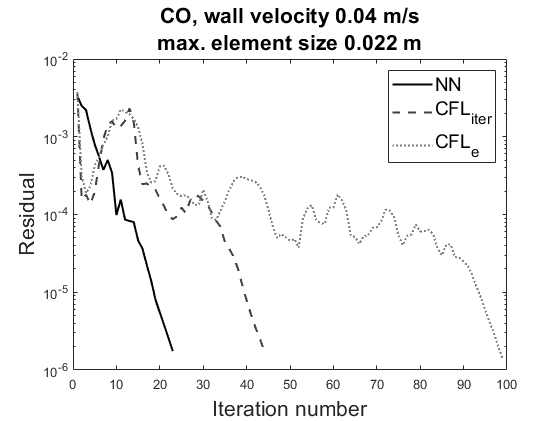}
	\caption{
%
	\AH{Convergence plot of a CO simulation with an ellipse with a-semiaxis 0.03\,m and b-semiaxis 0.05\,m, with center x = 0.1\,m and y = 0.4\,m for the local pseudo-time step under consideration: using the NN model as well as the strategies defined in~\cref{eq:CFLloc,eq:CFL2}.}
	\label{fig:convCO}}
\end{figure}

\begin{table}[t] \centering
\begin{tabular}{|p{1.9cm}|l|l|p{3.3cm}|p{2.8cm}|}
	\hline
	\textbf{geometry type} & \textbf{x (m)} & \textbf{y (m)}       & \textbf{inflow/wall \mbox{velocity} (\nicefrac{m}{s})} & \textbf{\mbox{max.\ element} \mbox{size (m)}} \\ \hline\hline
	\textbf{BO}            & 0.37            & 0.035\,:\,0.005\,:\,0.08 & 0.008                                      & 0.0126                          \\ \hline
	\textbf{BO}            & 0.3,\,1.1        & 0.04                 & 0.008                                      & 0.0126                          \\ \hline
	\textbf{CO}            & 0.1             & 0.4                  & 0.014\,:\,0.002\,:\,0.022                      & 0.01,\,0.03,\,0.05                \\ \hline
\end{tabular}\caption{Combinations of object center coordinates, velocities and mesh. Here, BO are B1 back-steps in \cref{tab:backstepgeom} with obstacles: a circle with radius 0.03\,m, an ellipse with a-semiaxis 0.03\,m and b-semiaxis 0.02\,m, and an ellipse with a-semiaxis 0.02\,m and b-semiaxis 0.03\,m. The Couette flow has geometry of C1 in \cref{tab:couettegeom} with obstacles: a circle with radius 0.04\,m, an ellipse with a-semiaxis 0.03\,m and b-semiaxis 0.02\,m, and an ellipse with a-semiaxis 0.02\,m and b-semiaxis 0.03\,m.}\label{tab:SimObst}
\end{table}

\begin{figure}[t]
	\begin{subfigure}{0.485\textwidth} \centering
		\input{Figures/Scatterplots/BO}
	\end{subfigure}
	\begin{subfigure}{0.485\textwidth} \centering
		\input{Figures/Scatterplots/CO}
	\end{subfigure}
	\caption{The results of the simulations with obstacles for each method; the simulations are ordered from smallest required number of nonlinear iterations to largest \AH{separately for each approach}. \AH{The horizontal lines indicate the average nonlinear iteration counts for the different approaches.}
		\label{fig:resultsO}
	}
\end{figure}
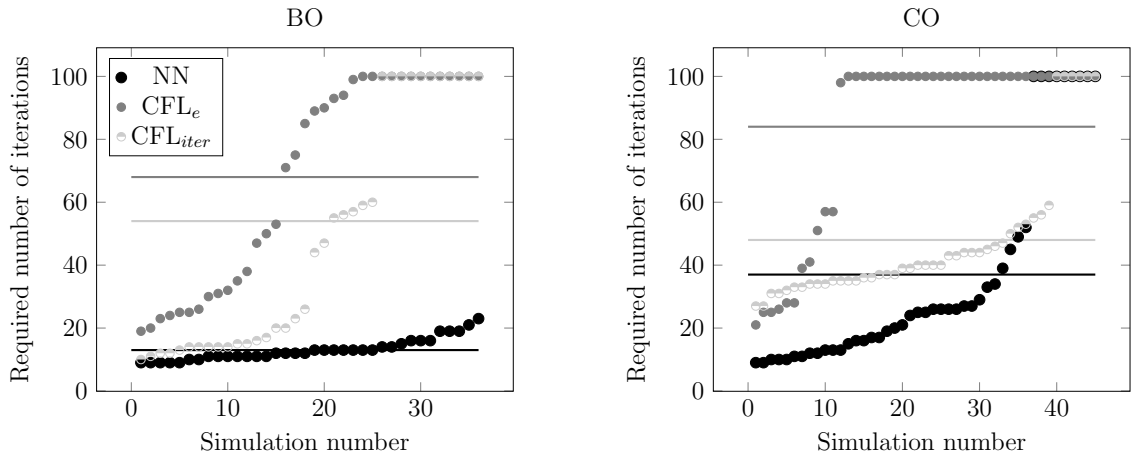

\subsection{Flow around an obstacle} \label{sec:results:obstacle}
As the final type of test examples, we consider Couette and back-step flow around an obstacle. In particular, we place an obstacle inside the computational domain of the back-step B1 and Couette C configurations; cf.~\Cref{tab:SimObst}\AH{, where} three different obstacle geometries are considered: a circle with radius of $0.03$\,m, an ellipse with a-semiaxis $0.03$\,m and b-semiaxis $0.02$\,m, and an ellipse with a-semiaxis $0.02$\,m and b-semiaxis $0.03$\,m. We denote the resulting configurations as back-step with obstacle (BO) and Couette cylinder with obstacle (CO); see~\cref{tab:SimObst} for more details on the configurations and~\cref{fig:couetteobst} for exemplary flow fields for the three different obstacle geometries in Couette flow. In total, we obtain $30$ BO and $45$ CO configurations.

The results of the simulations with obstacles are given in~\AZ{\Cref{fig:resultsO,fig:resultsCSO}}. Qualitatively, the NN approach compares to $\text{CFL}_{iter}$ and CFL$_e$ similarly as before. It performs better on $83.3$\,\% of the BO cases and on more than $73.3\,\%$ of the CO cases; it can also be observed that the CO cases are more difficult, as in more than $11$\,\% of the cases, none of the approaches converges. \AH{Also the average number of iteration is best for the NN approach for the BO and CO configurations.} For one of the considered simulations of CO with an elliptic obstacle, we again plot iterates and corresponding distributions of the local pseudo-time step in~\cref{fig:co}. In~\cref{fig:convCO}, we plot the residual history over the Newton iteration, showing the acceleration from the NN approach. 

\AH{Finally, also compared with the Newton methods in~\cref{sec:appendix}, the NN performs clearly best, that is, both in terms of robustness and in terms of acceleration of convergence.}

%% file: Figures/Scatterplots/BO.tex
\begin{tikzpicture}[scale=0.8]
\begin{axis}[
    title = BO,
    legend pos = north west,
    xlabel={Simulation number},
    ylabel={Required number of iterations}
    ]
    \addplot [
        scatter,
        only marks,
        point meta=explicit symbolic,
        scatter/classes={
            a={mark size=2.5pt,mark=*,black},
            b={mark size=2pt, mark=*,black!50},
            c={mark size=2pt, mark=myhalfcircle2,black!20}
        },
    ] table [meta=label] {
        x y label
        1 9 a
2 9 a
3 9 a
4 9 a
5 9 a
6 10 a
7 10 a
8 11 a
9 11 a
10 11 a
11 11 a
12 11 a
13 11 a
14 11 a
15 12 a
16 12 a
17 12 a
18 12 a
19 13 a
20 13 a
21 13 a
22 13 a
23 13 a
24 13 a
25 13 a
26 14 a
27 14 a
28 15 a
29 16 a
30 16 a
31 16 a
32 19 a
33 19 a
34 19 a
35 21 a
36 23 a
1 19 b
2 20 b
3 23 b
4 24 b
5 25 b
6 25 b
7 26 b
8 30 b
9 31 b
10 32 b
11 35 b
12 38 b
13 47 b
14 50 b
15 53 b
16 71 b
17 75 b
18 85 b
19 89 b
20 90 b
21 93 b
22 94 b
23 99 b
24 100 b
25 100 b
26 100 b
27 100 b
28 100 b
29 100 b
30 100 b
31 100 b
32 100 b
33 100 b
34 100 b
35 100 b
36 100 b
1 10 c
2 11 c
3 12 c
4 12 c
5 13 c
6 14 c
7 14 c
8 14 c
9 14 c
10 14 c
11 15 c
12 15 c
13 16 c
14 17 c
15 20 c
16 20 c
17 23 c
18 26 c
19 44 c
20 47 c
21 55 c
22 56 c
23 57 c
24 59 c
25 60 c
26 100 c
27 100 c
28 100 c
29 100 c
30 100 c
31 100 c
32 100 c
33 100 c
34 100 c
35 100 c
36 100 c
 };   
 \addplot[line width=1pt, color = black, 
] coordinates {(0,13) (36,13)};
\addplot[line width=1pt, color = black!50, 
] coordinates {(0,68) (36,68)};
\addplot[line width=1pt, color = black!20, 
] coordinates {(0,54) (36,54)};
\legend{NN,CFL$_e$,CFL$_{iter}$}
\end{axis}
\end{tikzpicture}

%% file: Figures/Scatterplots/CO.tex
\begin{tikzpicture}[scale=0.8]
\begin{axis}[
    title = CO,
    legend pos = north west,
    xlabel={Simulation number},
    ylabel={Required number of iterations}
    ]
    \addplot [
        scatter,
        only marks,
        point meta=explicit symbolic,
        scatter/classes={
            a={mark size=2.5pt,mark=*,black},
            b={mark size=2pt, mark=*,black!50},
            c={mark size=2pt, mark=myhalfcircle2,black!20}
        },
    ] table [meta=label] {
      x y label
      1 9 a
2 9 a
3 10 a
4 10 a
5 10 a
6 11 a
7 11 a
8 12 a
9 12 a
10 13 a
11 13 a
12 13 a
13 15 a
14 16 a
15 16 a
16 17 a
17 17 a
18 19 a
19 20 a
20 21 a
21 24 a
22 25 a
23 25 a
24 26 a
25 26 a
26 26 a
27 26 a
28 27 a
29 27 a
30 29 a
31 33 a
32 34 a
33 39 a
34 45 a
35 49 a
36 52 a
37 100 a
38 100 a
39 100 a
40 100 a
41 100 a
42 100 a
43 100 a
44 100 a
45 100 a
      1 21 b
2 25 b
3 25 b
4 26 b
5 28 b
6 28 b
7 39 b
8 41 b
9 51 b
10 57 b
11 57 b
12 98 b
13 100 b
14 100 b
15 100 b
16 100 b
17 100 b
18 100 b
19 100 b
20 100 b
21 100 b
22 100 b
23 100 b
24 100 b
25 100 b
26 100 b
27 100 b
28 100 b
29 100 b
30 100 b
31 100 b
32 100 b
33 100 b
34 100 b
35 100 b
36 100 b
37 100 b
38 100 b
39 100 b
40 100 b
41 100 b
42 100 b
43 100 b
44 100 b
45 100 b
1 27 c
2 27 c
3 31 c
4 31 c
5 32 c
6 33 c
7 33 c
8 34 c
9 34 c
10 34 c
11 35 c
12 35 c
13 35 c
14 35 c
15 36 c
16 36 c
17 37 c
18 37 c
19 37 c
20 39 c
21 39 c
22 40 c
23 40 c
24 40 c
25 40 c
26 43 c
27 43 c
28 44 c
29 44 c
30 44 c
31 45 c
32 46 c
33 47 c
34 50 c
35 52 c
36 53 c
37 55 c
38 56 c
39 59 c
40 100 c
41 100 c
42 100 c
43 100 c
44 100 c
45 100 c
}; 
\addplot[line width=1pt, color = black, 
] coordinates {(0,37) (45,37)};
\addplot[line width=1pt, color = black!50, 
] coordinates {(0,84) (45,84)};
\addplot[line width=1pt, color = black!20, 
] coordinates {(0,48) (45,48)};
\end{axis}
\end{tikzpicture}

%% file: Mainmatter/Conclusion.tex
\section{Conclusion}
\label{sec:conclusions}
\AH{We have introduced a novel} NN\AH{-based approach} to improve the pseudo-time stepping algorithm. Instead of using local pseudo-time steps based on predefined global CFL numbers, optimal local pseudo-time step predictions of the network are used to compute the local pseudo-time step. The network uses local information in order to make these local pseudo-time step predictions and to achieve generalizability. \AH{As a result, it can} accelerate the convergence for the simulations on which it has been trained as well as for simulations which it has not seen before. In all considered simulations, the network was able to perform better or equally well compared to the standard CFL number based pseudo-time stepping strategies in most cases\AH{; the same is true in comparison with the variants of Newton's method}. 

In future work, \AH{we would like to extend the method to} turbulent \AH{and three-dimensional flow cases}. \AH{Furthermore, we will} try to improve the method presented here, for example, by \AH{generating a more diverse training data set,} combining \AH{the heuristic strategies for choosing the pseudo-time step with our network approach} the default CFL number based strategies and the network, or removing redundan\AH{cies in the input features}. 

%% file: Mainmatter/Appendix.tex
\section{Comparison against variants of Newton's method}
\label{sec:appendix}

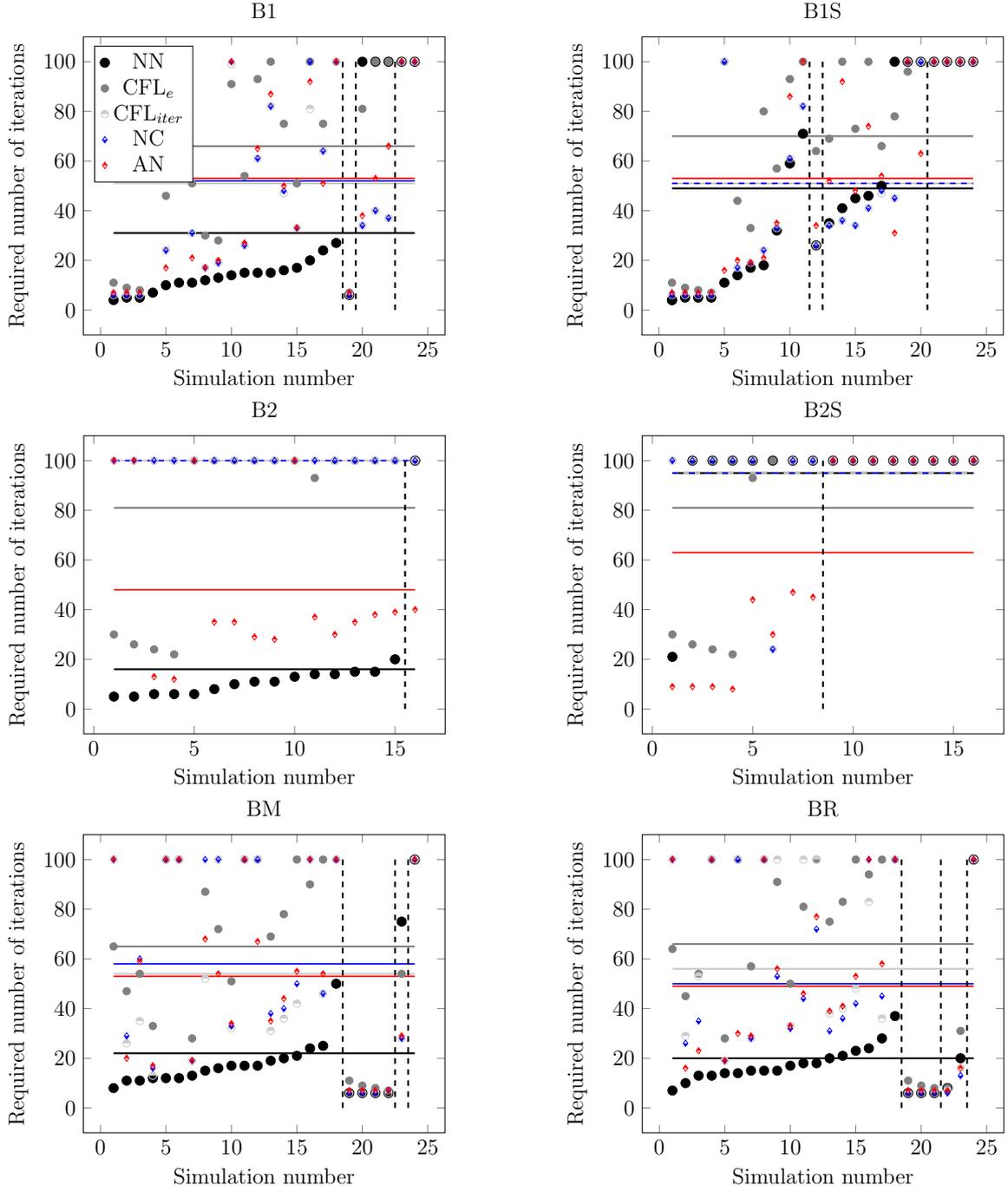
\begin{figure}[t]
	\begin{subfigure}{0.485\textwidth} \centering
		\input{Figures/Semi-ordered/B1}
	\end{subfigure}
	\begin{subfigure}{0.485\textwidth} \centering
		\input{Figures/Semi-ordered/B1S}
	\end{subfigure}
	\begin{subfigure}{0.485\textwidth} \centering
		\input{Figures/Semi-ordered/B2}
	\end{subfigure}
	\begin{subfigure}{0.485\textwidth} \centering
		\input{Figures/Semi-ordered/B2S}
	\end{subfigure}
	\begin{subfigure}{0.485\textwidth} \centering
		\input{Figures/Semi-ordered/BM}
	\end{subfigure}
	\hspace{9pt}
	\begin{subfigure}{0.485\textwidth} \centering
		\input{Figures/Semi-ordered/BR}
	\end{subfigure}
	\caption{\AH{Results for back-step flow configurations using the all approaches under consideration: using the NN, the strategies defined in~\cref{eq:CFLloc,eq:CFL2}, model as well as the NC and AN Newton methods; see~\cref{sec:appendix} for details on the organization of the plots.}
	\label{fig:resultsBSO}}
\end{figure}

\begin{figure}[t]
	\begin{subfigure}{0.485\textwidth} \centering
		\input{Figures/Semi-ordered/C}
	\end{subfigure}
	\begin{subfigure}{0.485\textwidth} \centering
		\input{Figures/Semi-ordered/CS}
	\end{subfigure}
	\begin{subfigure}{0.485\textwidth} \centering
		\input{Figures/Semi-ordered/BO}
	\end{subfigure}
	\hspace{4pt}
	\begin{subfigure}{0.485\textwidth} \centering
		\input{Figures/Semi-ordered/CO}
	\end{subfigure}
	\caption{\AH{Results for Couette flow and flow around an obstacle configurations using the all approaches under consideration: using the NN, the strategies defined in~\cref{eq:CFLloc,eq:CFL2}, model as well as the NC and AN Newton methods; see~\cref{sec:appendix} for details on the organization of the plots. The horizontal lines indicate the average nonlinear iteration counts for the different approaches.}
	\label{fig:resultsCSO}}
\end{figure}
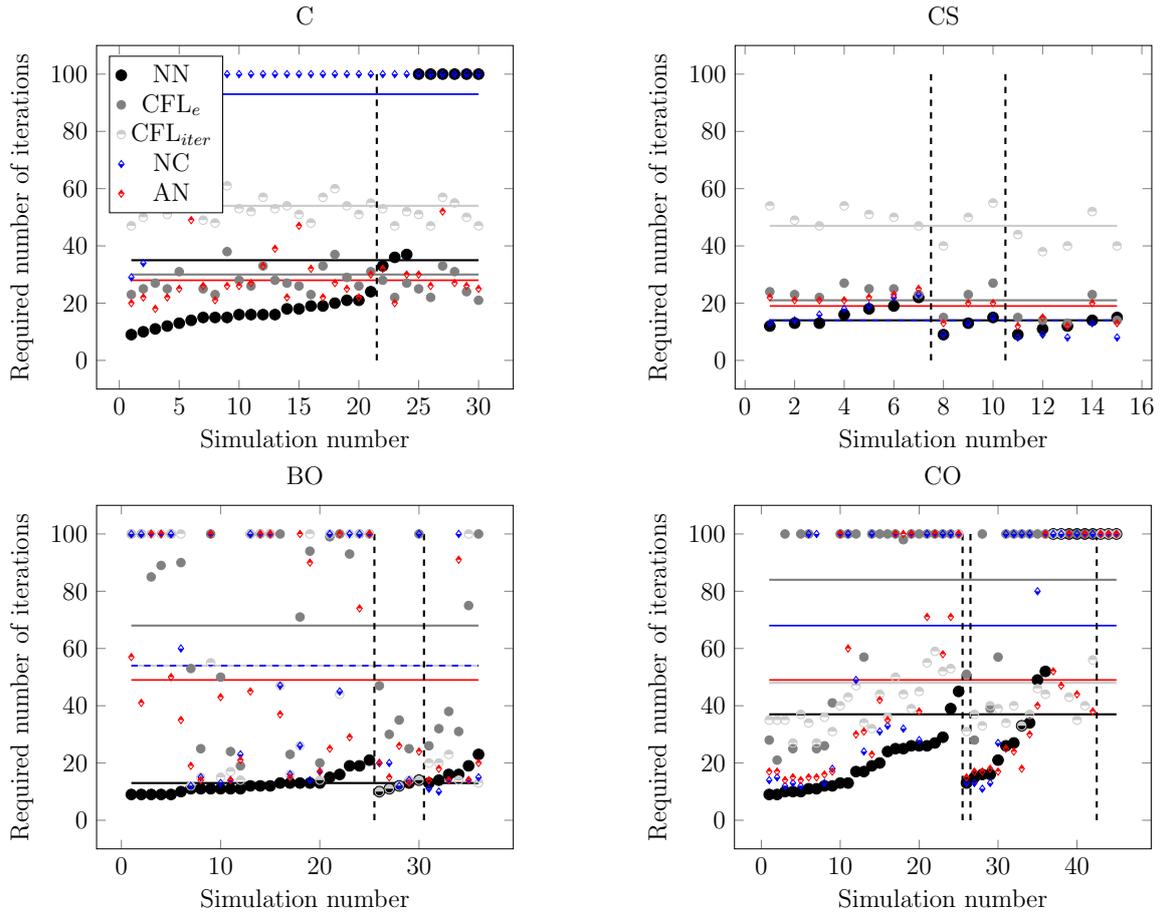

In~\Cref{fig:resultsBSO,fig:resultsCSO}, we show results for all configurations \AH{considered, comparing our novel NN approach with four nonlinear solvers implemented in COMSOL, that is, the two default CFL number strategies defined in~\cref{eq:CFLloc,eq:CFL2} as well as two variants of Newton's method. As variants of Newton's method, we employ the standard Newton method with damping factor of one and an adaptive choice of the damping factor, depending on the intermediate iterates. We denote the two approaches as Newton constant (NC) and automatic Newton (AN), respectively; cf.~\cite[p.~1627]{CMRM}. 
Different from~\Cref{fig:resultsB,fig:resultsC,fig:resultsO}, the ordering of the simulations for the different methods is the same for each approach. The simulations are first organized into different blocks, and within each block the ordering is based on the number of nonlinear iterations needed for the NN approach. The blocks are as follows: before the first vertical dashed line, the NN performed best; between the first and second vertical dashed line, the NN performed best with the same number of iterations as at least one other method; between the second and third vertical dashed line, the network performed worse than at least one other method; after the third vertical line, none of the methods converged.}

\AH{The results are mostly in alignment with the previous results as the NC and AN methods often perform similarly to one of the classical CFL strategies considered. In particular, when inspecting the average iteration count, both NC and AN perform quite similarly to CFL$_{iter}$ for B1, B1S, BM, BR, CS, and BO configurations; in all those cases the NN approach is clearly competitive. For the CO and B2 cases, AN performs better than NC, while the NN approach performs best. Only for the B2S configuration, where the NN approach generally performed worst, AC seems to perform clearly best; however, all approaches struggle for the B2S cases. Finally, for the C configurations, all methods perform similarly well, except for NC, which fails in most of the cases.}

\AH{In summary, adding the NC and AN approaches to the picture does not change the qualitative comparison much; only for the B2S cases, AN performs clearly better than the NN approach. In future work, this should be taken into account for further improving the NN model.}


%% file: Figures/Semi-ordered/B1.tex
\begin{tikzpicture}[scale=0.8]
\begin{axis}[
    title = B1,
    legend pos = north west,
    xlabel={Simulation number},
    ylabel={Required number of iterations}
    ]
    \addplot [
        scatter,
        only marks,
        point meta=explicit symbolic,
        scatter/classes={
            a={mark size=2.5pt,mark=*,black},
            b={mark size=2pt, mark=*,black!50},
            c={mark size=2pt, mark=myhalfcircle2,black!20},
            d={mark size=1.7pt, mark=mymark1,blue},
            e={mark size=1.7pt, mark=mymark2,red}
        },
    ] table [meta=label] {
        x      y        label
1	4	a
2	5	a
3	5	a
4	7	a
5	10	a
6	11	a
7	11	a
8	12	a
9	13	a
10	14	a
11	15	a
12	15	a
13	15	a
14	16	a
15	17	a
16	20	a
17	24	a
18	27	a
19	6	a
20	100	a
21	100	a
22	100	a
23	100	a
24	100	a
1	11	b
2	9	b
3	8	b
4	65	b
5	46	b
6	100	b
7	51	b
8	30	b
9	28	b
10	91	b
11	54	b
12	93	b
13	100	b
14	75	b
15	51	b
16	100	b
17	75	b
18	100	b
19	7	b
20	81	b
21	100	b
22	100	b
23	100	b
24	100	b
1	6	c
2	6	c
3	6	c
4	100	c
5	24	c
6	100	c
7	31	c
8	17	c
9	19	c
10	99	c
11	26	c
12	61	c
13	82	c
14	47	c
15	33	c
16	81	c
17	64	c
18	100	c
19	6	c
20	34	c
21	40	c
22	37	c
23	100	c
24	100	c
1	6	d
2	6	d
3	6	d
4	100	d
5	24	d
6	100	d
7	31	d
8	17	d
9	19	d
10	100	d
11	26	d
12	61	d
13	82	d
14	48	d
15	33	d
16	100	d
17	64	d
18	100	d
19	6	d
20	34	d
21	40	d
22	37	d
23	100	d
24	100	d
1	7	e
2	7	e
3	7	e
4	100	e
5	17	e
6	100	e
7	21	e
8	17	e
9	20	e
10	100	e
11	27	e
12	65	e
13	87	e
14	50	e
15	33	e
16	92	e
17	51	e
18	100	e
19	7	e
20	38	e
21	53	e
22	66	e
23	100	e
24	100	e
    };  
    
\addplot[line width=1pt, color = black, 
] coordinates {(1,31) (24,31)};
\addplot[line width=1pt, color = black!50, 
] coordinates {(1,66) (24,66)};
\addplot[line width=1pt, color = black!20, 
] coordinates {(1,51) (24,51)};
\legend{NN,CFL$_r$,CFL$_{iter}$,NC,AN}
\addplot[line width=0.8pt, color = blue, 
] coordinates {(1,52) (24,52)};
\legend{NN,CFL$_r$,CFL$_{iter}$,NC,AN}
\addplot[line width=0.8pt, color = red, 
] coordinates {(1,53) (24,53)};
\legend{NN,CFL$_e$,CFL$_{iter}$,NC,AN}

\addplot[line width=1pt, color = black, 
style = dashed
] coordinates {(18.5,0) (18.5,100)};
\addplot[line width=1pt, color = black, 
style = dashed
] coordinates {(19.5,0) (19.5,100)};
\addplot[line width=1pt, color = black, 
style = dashed
] coordinates {(22.5,0) (22.5,100)};

\end{axis}
\end{tikzpicture}

%% file: Figures/Semi-ordered/B1S.tex
\begin{tikzpicture}[scale=0.8]
\begin{axis}[
    title = B1S,
    legend pos = north west,
    xlabel={Simulation number},
    ylabel={Required number of iterations}
    ]
    \addplot [
        scatter,
        only marks,
        point meta=explicit symbolic,
        scatter/classes={
            a={mark size=2.5pt,mark=*,black},
            b={mark size=2pt, mark=*,black!50},
            c={mark size=2pt, mark=myhalfcircle2,black!20},
            d={mark size=1.7pt, mark=mymark1,blue},
            e={mark size=1.7pt, mark=mymark2,red}
        },
    ] table [meta=label] {
        x y label
        1 4 a
1	4	a
2	5	a
3	5	a
4	5	a
5	11	a
6	14	a
7	17	a
8	18	a
9	32	a
10	59	a
11	71	a
12	26	a
13	35	a
14	41	a
15	45	a
16	46	a
17	50	a
18	100	a
19	100	a
20	100	a
21	100	a
22	100	a
23	100	a
24	100	a
1	11	b
2	9	b
3	8	b
4	7	b
5	100	b
6	44	b
7	33	b
8	80	b
9	57	b
10	93	b
11	100	b
12	64	b
13	69	b
14	100	b
15	73	b
16	100	b
17	66	b
18	78	b
19	96	b
20	100	b
21	100	b
22	100	b
23	100	b
24	100	b
1	6	c
2	6	c
3	6	c
4	6	c
5	100	c
6	17	c
7	19	c
8	24	c
9	33	c
10	61	c
11	82	c
12	26	c
13	34	c
14	36	c
15	34	c
16	41	c
17	47	c
18	45	c
19	100	c
20	100	c
21	100	c
22	100	c
23	100	c
24	100	c
1	6	d
2	6	d
3	6	d
4	6	d
5	100	d
6	17	d
7	19	d
8	24	d
9	33	d
10	61	d
11	82	d
12	26	d
13	34	d
14	36	d
15	34	d
16	41	d
17	48	d
18	45	d
19	100	d
20	100	d
21	100	d
22	100	d
23	100	d
24	100	d
1	7	e
2	7	e
3	7	e
4	7	e
5	16	e
6	20	e
7	19	e
8	21	e
9	35	e
10	86	e
11	100	e
12	34	e
13	52	e
14	92	e
15	48	e
16	74	e
17	54	e
18	31	e
19	100	e
20	63	e
21	100	e
22	100	e
23	100	e
24	100	e

};

\addplot[line width=1pt, color = black 
] coordinates {(1,49) (24,49)};
\addplot[line width=1pt, color = black!50 
] coordinates {(1,70) (24,70)};
\addplot[line width=1pt, color = black!20 
] coordinates {(1,51) (24,51)};
\addplot[line width=0.8pt, color = blue, 
style = dashed
] coordinates {(1,51) (24,51)};
\addplot[line width=0.8pt, color = red, 
] coordinates {(1,53) (24,53)};

\addplot[line width=1pt, color = black, 
style = dashed
] coordinates {(11.5,0) (11.5,100)};
\addplot[line width=1pt, color = black, 
style = dashed
] coordinates {(12.5,0) (12.5,100)};
\addplot[line width=1pt, color = black,
style = dashed
] coordinates {(20.5,0) (20.5,100)};

\end{axis}
\end{tikzpicture}

%% file: Figures/Semi-ordered/B2.tex
\begin{tikzpicture}[scale=0.8]
\begin{axis}[
    title = B2,
    xlabel={Simulation number},
    ylabel={Required number of iterations}
    ]
    \addplot [
        scatter,
        only marks,
        point meta=explicit symbolic,
        scatter/classes={
            a={mark size=2.5pt,mark=*,black},
            b={mark size=2pt, mark=*,black!50},
            c={mark size=2pt, mark=myhalfcircle2,black!20},
            d={mark size=1.7pt, mark=mymark1,blue},
            e={mark size=1.7pt, mark=mymark2,red}
        },
    ] table [meta=label] {
        x y label
1	5	a
2	5	a
3	6	a
4	6	a
5	6	a
6	8	a
7	10	a
8	11	a
9	11	a
10	13	a
11	14	a
12	14	a
13	15	a
14	15	a
15	20	a
16	100	a
1	30	b
2	26	b
3	24	b
4	22	b
5	100	b
6	100	b
7	100	b
8	100	b
9	100	b
10	100	b
11	93	b
12	100	b
13	100	b
14	100	b
15	100	b
16	100	b
1	100	c
2	100	c
3	100	c
4	100	c
5	100	c
6	100	c
7	100	c
8	100	c
9	100	c
10	100	c
11	100	c
12	100	c
13	100	c
14	100	c
15	100	c
16	100	c
1	100	d
2	100	d
3	100	d
4	100	d
5	100	d
6	100	d
7	100	d
8	100	d
9	100	d
10	100	d
11	100	d
12	100	d
13	100	d
14	100	d
15	100	d
16	100	d
1	100	e
2	100	e
3	13	e
4	12	e
5	100	e
6	35	e
7	35	e
8	29	e
9	28	e
10	100	e
11	37	e
12	30	e
13	35	e
14	38	e
15	39	e
16	40	e
};    
    \addplot[line width=1pt, color = black, 
] coordinates {(1,16) (16,16)};
\addplot[line width=1pt, color = black!50, 
] coordinates {(1,81) (16,81)};
\addplot[line width=1pt, color = black!20, 
] coordinates {(1,100) (16,100)};
\addplot[line width=0.8pt, color = blue, 
style = dashed
] coordinates {(1,100) (16,100)};
\addplot[line width=0.8pt, color = red, 
] coordinates {(1,48) (16,48)};

\addplot[line width=1pt, color = black, 
style = dashed
] coordinates {(15.5,0) (15.5,100)};

\end{axis}
\end{tikzpicture}

%% file: Figures/Semi-ordered/B2S.tex
\begin{tikzpicture}[scale=0.8]
\begin{axis}[
    title = B2S,
    legend pos = south east,
    xlabel={Simulation number},
    ylabel={Required number of iterations}
    ]
    \addplot [
        scatter,
        only marks,
        point meta=explicit symbolic,
        scatter/classes={
            a={mark size=2.5pt,mark=*,black},
            b={mark size=2pt, mark=*,black!50},
            c={mark size=2pt, mark=myhalfcircle2,black!20},
            d={mark size=1.7pt, mark=mymark1,blue},
            e={mark size=1.7pt, mark=mymark2,red}
        },
    ] table [meta=label] {
        x y label
        1	21	a
2	100	a
3	100	a
4	100	a
5	100	a
6	100	a
7	100	a
8	100	a
9	100	a
10	100	a
11	100	a
12	100	a
13	100	a
14	100	a
15	100	a
16	100	a
1	30	b
2	26	b
3	24	b
4	22	b
5	93	b
6	100	b
7	100	b
8	100	b
9	100	b
10	100	b
11	100	b
12	100	b
13	100	b
14	100	b
15	100	b
16	100	b
1	100	c
2	100	c
3	100	c
4	100	c
5	100	c
6	24	c
7	100	c
8	100	c
9	100	c
10	100	c
11	100	c
12	100	c
13	100	c
14	100	c
15	100	c
16	100	c
1	100	d
2	100	d
3	100	d
4	100	d
5	100	d
6	24	d
7	100	d
8	100	d
9	100	d
10	100	d
11	100	d
12	100	d
13	100	d
14	100	d
15	100	d
16	100	d
1	9	e
2	9	e
3	9	e
4	8	e
5	44	e
6	30	e
7	47	e
8	45	e
9	100	e
10	100	e
11	100	e
12	100	e
13	100	e
14	100	e
15	100	e
16	100	e
};    
    \addplot[line width=1pt, color = black, 
] coordinates {(1,95) (16,95)};
\addplot[line width=1pt, color = black!50, 
] coordinates {(1,81) (16,81)};
\addplot[line width=1pt, color = black!20, 
style = dashed
] coordinates {(1,95) (16,95)};
\addplot[line width=0.8pt, color = blue, 
style = loosely dashed
] coordinates {(1,95) (16,95)};
\addplot[line width=0.8pt, color = red, 
] coordinates {(1,63) (16,63)};

\addplot[line width=1pt, color = black, 
style = dashed
] coordinates {(8.5,0) (8.5,100)};

\end{axis}
\end{tikzpicture}

%% file: Figures/Semi-ordered/BM.tex
\begin{tikzpicture}[scale=0.8]
\begin{axis}[
    title = BM,
    legend pos = north west,
    xlabel={Simulation number},
    ylabel={Required number of iterations}
    ]
    \addplot [
        scatter,
        only marks,
        point meta=explicit symbolic,
        scatter/classes={
            a={mark size=2.5pt,mark=*,black},
            b={mark size=2pt, mark=*,black!50},
            c={mark size=2pt, mark=myhalfcircle2,black!20},
            d={mark size=1.7pt, mark=mymark1,blue},
            e={mark size=1.7pt, mark=mymark2,red}
        },
    ] table [meta=label] {
        x y label
1	8	a
2	11	a
3	11	a
4	12	a
5	12	a
6	12	a
7	13	a
8	15	a
9	16	a
10	17	a
11	17	a
12	17	a
13	19	a
14	20	a
15	21	a
16	24	a
17	25	a
18	50	a
19	6	a
20	6	a
21	6	a
22	6	a
23	75	a
24	100	a
1	65	b
2	47	b
3	54	b
4	33	b
5	100	b
6	100	b
7	28	b
8	87	b
9	72	b
10	51	b
11	100	b
12	100	b
13	69	b
14	78	b
15	100	b
16	90	b
17	100	b
18	100	b
19	11	b
20	9	b
21	8	b
22	7	b
23	54	b
24	100	b
1	100	c
2	26	c
3	35	c
4	13	c
5	100	c
6	100	c
7	19	c
8	52	c
9	100	c
10	32	c
11	100	c
12	100	c
13	31	c
14	36	c
15	42	c
16	100	c
17	46	c
18	100	c
19	6	c
20	6	c
21	6	c
22	6	c
23	28	c
24	100	c
1	100	d
2	29	d
3	60	d
4	16	d
5	100	d
6	100	d
7	19	d
8	100	d
9	100	d
10	33	d
11	100	d
12	100	d
13	38	d
14	40	d
15	50	d
16	100	d
17	46	d
18	100	d
19	6	d
20	6	d
21	6	d
22	7	d
23	28	d
24	100	d
1	100	e
2	20	e
3	59	e
4	17	e
5	100	e
6	100	e
7	19	e
8	68	e
9	54	e
10	34	e
11	100	e
12	67	e
13	35	e
14	44	e
15	55	e
16	100	e
17	54	e
18	100	e
19	7	e
20	7	e
21	7	e
22	7	e
23	29	e
24	100	e
};    
 \addplot[line width=1pt, color = black, 
] coordinates {(1,22) (24,22)};
\addplot[line width=1pt, color = black!50, 
] coordinates {(1,65) (24,65)};
\addplot[line width=1pt, color = black!20, 
] coordinates {(1,54) (24,54)};
\addplot[line width=0.8pt, color = blue, 
] coordinates {(1,58) (24,58)};
\addplot[line width=0.8pt, color = red, 
] coordinates {(1,53) (24,53)};

\addplot[line width=1pt, color = black, 
style = dashed
] coordinates {(18.5,0) (18.5,100)};
\addplot[line width=1pt, color = black, 
style = dashed
] coordinates {(22.5,0) (22.5,100)};
\addplot[line width=1pt, color = black, 
style = dashed
] coordinates {(23.5,0) (23.5,100)};
\end{axis}
\end{tikzpicture}

%% file: Figures/Semi-ordered/BR.tex
\begin{tikzpicture}[scale=0.8]
\begin{axis}[
    title = BR,
    legend pos = north west,
    xlabel={Simulation number},
    ylabel={Required number of iterations}
    ]
    \addplot [
        scatter,
        only marks,
        point meta=explicit symbolic,
        scatter/classes={
            a={mark size=2.5pt,mark=*,black},
            b={mark size=2pt, mark=*,black!50},
            c={mark size=2pt, mark=myhalfcircle2,black!20},
            d={mark size=1.7pt, mark=mymark1,blue},
            e={mark size=1.7pt, mark=mymark2,red}
        },
    ] table [meta=label] {
        x y label
1	7	a
2	10	a
3	13	a
4	13	a
5	14	a
6	14	a
7	15	a
8	15	a
9	15	a
10	17	a
11	18	a
12	18	a
13	20	a
14	21	a
15	23	a
16	24	a
17	28	a
18	37	a
19	6	a
20	6	a
21	6	a
22	8	a
23	20	a
24	100	a
1	64	b
2	45	b
3	54	b
4	100	b
5	28	b
6	100	b
7	57	b
8	100	b
9	91	b
10	50	b
11	81	b
12	100	b
13	75	b
14	83	b
15	100	b
16	94	b
17	100	b
18	100	b
19	11	b
20	9	b
21	8	b
22	8	b
23	31	b
24	100	b
1	100	c
2	29	c
3	53	c
4	100	c
5	19	c
6	100	c
7	28	c
8	100	c
9	100	c
10	33	c
11	100	c
12	100	c
13	38	c
14	40	c
15	48	c
16	83	c
17	36	c
18	100	c
19	6	c
20	6	c
21	6	c
22	7	c
23	16	c
24	100	c
1	100	d
2	26	d
3	35	d
4	100	d
5	19	d
6	100	d
7	28	d
8	100	d
9	53	d
10	32	d
11	44	d
12	72	d
13	31	d
14	36	d
15	42	d
16	100	d
17	45	d
18	100	d
19	6	d
20	6	d
21	6	d
22	6	d
23	13	d
24	100	d
1	100	e
2	16	e
3	23	e
4	100	e
5	19	e
6	30	e
7	29	e
8	100	e
9	56	e
10	33	e
11	46	e
12	77	e
13	39	e
14	41	e
15	53	e
16	100	e
17	58	e
18	100	e
19	7	e
20	7	e
21	7	e
22	7	e
23	16	e
24	100	e
};    
 \addplot[line width=1pt, color = black, 
] coordinates {(1,20) (24,20)};
\addplot[line width=1pt, color = black!50, 
] coordinates {(1,66) (24,66)};
\addplot[line width=1pt, color = black!20, 
] coordinates {(1,56) (24,56)};
\addplot[line width=0.8pt, color = blue, 
] coordinates {(1,50) (24,50)};
\addplot[line width=0.8pt, color = red, 
] coordinates {(1,49) (24,49)};

\addplot[line width=1pt, color = black, 
style = dashed
] coordinates {(18.5,0) (18.5,100)};
\addplot[line width=1pt, color = black, 
style = dashed
] coordinates {(21.5,0) (21.5,100)};
\addplot[line width=1pt, color = black, 
style = dashed
] coordinates {(23.5,0) (23.5,100)};
\end{axis}
\end{tikzpicture}

%% file: Figures/Semi-ordered/C.tex
\begin{tikzpicture}[scale=0.8]
\begin{axis}[
    title = C,
    legend pos = north west,
    xlabel={Simulation number},
    ylabel={Required number of iterations}
    ]
    \addplot [
        scatter,
        only marks,
        point meta=explicit symbolic,
        scatter/classes={
            a={mark size=2.5pt,mark=*,black},
            b={mark size=2pt, mark=*,black!50},
            c={mark size=2pt, mark=myhalfcircle2,black!20},
            d={mark size=1.7pt, mark=mymark1,blue},
            e={mark size=1.7pt, mark=mymark2,red}
        },
    ] table [meta=label] {
        x y label
1	9	a
2	10	a
3	11	a
4	12	a
5	13	a
6	14	a
7	15	a
8	15	a
9	15	a
10	16	a
11	16	a
12	16	a
13	16	a
14	18	a
15	18	a
16	19	a
17	19	a
18	20	a
19	21	a
20	21	a
21	24	a
22	33	a
23	36	a
24	37	a
25	100	a
26	100	a
27	100	a
28	100	a
29	100	a
30	100	a
1	23	b
2	25	b
3	27	b
4	25	b
5	31	b
6	100	b
7	25	b
8	23	b
9	38	b
10	28	b
11	26	b
12	33	b
13	28	b
14	27	b
15	26	b
16	23	b
17	33	b
18	37	b
19	29	b
20	26	b
21	31	b
22	28	b
23	22	b
24	27	b
25	25	b
26	22	b
27	33	b
28	31	b
29	24	b
30	21	b
1	47	c
2	50	c
3	55	c
4	51	c
5	55	c
6	100	c
7	49	c
8	48	c
9	61	c
10	53	c
11	52	c
12	57	c
13	53	c
14	54	c
15	51	c
16	48	c
17	57	c
18	60	c
19	54	c
20	51	c
21	55	c
22	53	c
23	47	c
24	52	c
25	51	c
26	47	c
27	57	c
28	55	c
29	50	c
30	47	c
1	29	d
2	34	d
3	100	d
4	56	d
5	100	d
6	100	d
7	100	d
8	63	d
9	100	d
10	100	d
11	100	d
12	100	d
13	100	d
14	100	d
15	100	d
16	100	d
17	100	d
18	100	d
19	100	d
20	100	d
21	100	d
22	100	d
23	100	d
24	100	d
25	100	d
26	100	d
27	100	d
28	100	d
29	100	d
30	100	d
1	20	e
2	22	e
3	18	e
4	22	e
5	25	e
6	49	e
7	26	e
8	21	e
9	26	e
10	26	e
11	27	e
12	33	e
13	39	e
14	22	e
15	47	e
16	32	e
17	22	e
18	27	e
19	25	e
20	22	e
21	30	e
22	32	e
23	20	e
24	30	e
25	30	e
26	26	e
27	52	e
28	27	e
29	26	e
30	25	e
};  
  \addplot[line width=1pt, color = black, 
] coordinates {(1,35) (30,35)};
\addplot[line width=1pt, color = black!50, 
] coordinates {(1,30) (30,30)};
\addplot[line width=1pt, color = black!20, 
] coordinates {(1,54) (30,54)};
\legend{NN,CFL$_e$,CFL$_{iter}$,NC,AN}
\addplot[line width=0.8pt, color = blue, 
] coordinates {(1,93) (30,93)};
\addplot[line width=0.8pt, color = red, 
] coordinates {(1,28) (30,28)};

\addplot[line width=1pt, color = black, 
style = dashed
] coordinates {(21.5,0) (21.5,100)};

\end{axis}
\end{tikzpicture}

%% file: Figures/Semi-ordered/CS.tex
\begin{tikzpicture}[scale=0.8]
\begin{axis}[
    title = CS,
    legend pos = north west,
    xlabel={Simulation number},
    ylabel={Required number of iterations}
    ]
    \addplot [
        scatter,
        only marks,
        point meta=explicit symbolic,
        scatter/classes={
            a={mark size=2.5pt,mark=*,black},
            b={mark size=2pt, mark=*,black!50},
            c={mark size=2pt, mark=myhalfcircle2,black!20},
            d={mark size=1.7pt, mark=mymark1,blue},
            e={mark size=1.7pt, mark=mymark2,red}
        },
    ] table [meta=label] {
       x y label
1	12	a
2	13	a
3	13	a
4	16	a
5	18	a
6	19	a
7	22	a
8	9	a
9	13	a
10	15	a
11	9	a
12	11	a
13	12	a
14	14	a
15	15	a
1	24	b
2	23	b
3	22	b
4	27	b
5	25	b
6	25	b
7	23	b
8	15	b
9	23	b
10	27	b
11	15	b
12	14	b
13	13	b
14	23	b
15	14	b
1	54	c
2	49	c
3	47	c
4	54	c
5	51	c
6	50	c
7	47	c
8	40	c
9	50	c
10	55	c
11	44	c
12	38	c
13	40	c
14	52	c
15	40	c
1	13	d
2	14	d
3	16	d
4	18	d
5	19	d
6	22	d
7	23	d
8	9	d
9	13	d
10	15	d
11	8	d
12	9	d
13	8	d
14	13	d
15	8	d
1	22	e
2	21	e
3	21	e
4	21	e
5	22	e
6	23	e
7	25	e
8	13	e
9	20	e
10	20	e
11	12	e
12	15	e
13	12	e
14	20	e
15	13	e
};    
  \addplot[line width=1pt, color = black, 
] coordinates {(1,14) (15,14)};
\addplot[line width=1pt, color = black!50, 
] coordinates {(1,21) (15,21)};
\addplot[line width=1pt, color = black!20, 
] coordinates {(1,47) (15,47)};
\addplot[line width=0.8pt, color = blue, 
style = dashed
] coordinates {(1,14) (15,14)};
\addplot[line width=0.8pt, color = red, 
] coordinates {(1,19) (15,19)};

\addplot[line width=1pt, color = black, 
style = dashed
] coordinates {(7.5,0) (7.5,100)};
\addplot[line width=1pt, color = black, 
style = dashed
] coordinates {(10.5,0) (10.5,100)};

\end{axis}
\end{tikzpicture}

%% file: Figures/Semi-ordered/BO.tex
\begin{tikzpicture}[scale=0.8]
\begin{axis}[
    title = BO,
    legend pos = north west,
    xlabel={Simulation number},
    ylabel={Required number of iterations}
    ]
    \addplot [
        scatter,
        only marks,
        point meta=explicit symbolic,
        scatter/classes={
            a={mark size=2.5pt,mark=*,black},
            b={mark size=2pt, mark=*,black!50},
            c={mark size=2pt, mark=myhalfcircle2,black!20},
            d={mark size=1.7pt, mark=mymark1,blue},
            e={mark size=1.7pt, mark=mymark2,red}
        },
    ] table [meta=label] {
        x y label
1	9	a
2	9	a
3	9	a
4	9	a
5	9	a
6	10	a
7	11	a
8	11	a
9	11	a
10	11	a
11	11	a
12	11	a
13	12	a
14	12	a
15	12	a
16	13	a
17	13	a
18	13	a
19	13	a
20	13	a
21	15	a
22	16	a
23	19	a
24	19	a
25	21	a
26	10	a
27	11	a
28	12	a
29	13	a
30	14	a
31	13	a
32	14	a
33	16	a
34	16	a
35	19	a
36	23	a
1	100	b
2	100	b
3	85	b
4	89	b
5	100	b
6	90	b
7	53	b
8	25	b
9	100	b
10	50	b
11	24	b
12	19	b
13	100	b
14	100	b
15	100	b
16	100	b
17	23	b
18	71	b
19	94	b
20	20	b
21	99	b
22	100	b
23	93	b
24	100	b
25	100	b
26	47	b
27	30	b
28	35	b
29	25	b
30	100	b
31	26	b
32	32	b
33	38	b
34	31	b
35	75	b
36	100	b
1	100	c
2	100	c
3	100	c
4	100	c
5	100	c
6	100	c
7	12	c
8	15	c
9	55	c
10	15	c
11	17	c
12	14	c
13	100	c
14	100	c
15	100	c
16	47	c
17	16	c
18	26	c
19	100	c
20	14	c
21	100	c
22	44	c
23	100	c
24	100	c
25	100	c
26	10	c
27	11	c
28	12	c
29	14	c
30	14	c
31	20	c
32	20	c
33	23	c
34	14	c
35	100	c
36	13	c
1	100	d
2	100	d
3	100	d
4	100	d
5	100	d
6	60	d
7	12	d
8	15	d
9	100	d
10	13	d
11	14	d
12	23	d
13	100	d
14	100	d
15	100	d
16	47	d
17	16	d
18	26	d
19	14	d
20	17	d
21	100	d
22	45	d
23	100	d
24	100	d
25	100	d
26	20	d
27	20	d
28	12	d
29	14	d
30	100	d
31	11	d
32	10	d
33	14	d
34	100	d
35	14	d
36	15	d
1	57	e
2	41	e
3	100	e
4	100	e
5	50	e
6	35	e
7	19	e
8	14	e
9	100	e
10	43	e
11	14	e
12	21	e
13	45	e
14	100	e
15	100	e
16	37	e
17	15	e
18	100	e
19	90	e
20	17	e
21	25	e
22	100	e
23	29	e
24	74	e
25	100	e
26	20	e
27	15	e
28	26	e
29	13	e
30	24	e
31	14	e
32	18	e
33	14	e
34	91	e
35	14	e
36	20	e
};    
 \addplot[line width=1pt, color = black, 
] coordinates {(1,13) (36,13)};
\addplot[line width=1pt, color = black!50, 
] coordinates {(1,68) (36,68)};
\addplot[line width=1pt, color = black!20, 
] coordinates {(1,54) (36,54)};
\addplot[line width=0.8pt, color = blue, 
style = dashed
] coordinates {(1,54) (36,54)};
\addplot[line width=0.8pt, color = red, 
] coordinates {(1,49) (36,49)};

\addplot[line width=1pt, color = black, 
style = dashed
] coordinates {(25.5,0) (25.5,100)};
\addplot[line width=1pt, color = black, 
style = dashed
] coordinates {(30.5,0) (30.5,100)};
\end{axis}
\end{tikzpicture}

%% file: Figures/Semi-ordered/CO.tex
\begin{tikzpicture}[scale=0.8]
\begin{axis}[
    title = CO,
    legend pos = north west,
    xlabel={Simulation number},
    ylabel={Required number of iterations}
    ]
    \addplot [
        scatter,
        only marks,
        point meta=explicit symbolic,
        scatter/classes={
            a={mark size=2.5pt,mark=*,black},
            b={mark size=2pt, mark=*,black!50},
            c={mark size=2pt, mark=myhalfcircle2,black!20},
            d={mark size=1.7pt, mark=mymark1,blue},
            e={mark size=1.7pt, mark=mymark2,red}
        },
    ] table [meta=label] {
x y label
1	9	a
2	9	a
3	10	a
4	10	a
5	10	a
6	11	a
7	11	a
8	12	a
9	12	a
10	13	a
11	13	a
12	17	a
13	17	a
14	19	a
15	20	a
16	24	a
17	25	a
18	25	a
19	26	a
20	26	a
21	26	a
22	27	a
23	29	a
24	39	a
25	45	a
26	13	a
27	15	a
28	16	a
29	16	a
30	21	a
31	26	a
32	27	a
33	33	a
34	34	a
35	49	a
36	52	a
37	100	a
38	100	a
39	100	a
40	100	a
41	100	a
42	100	a
43	100	a
44	100	a
45	100	a
1	28	b
2	21	b
3	100	b
4	25	b
5	100	b
6	100	b
7	25	b
8	26	b
9	41	b
10	100	b
11	100	b
12	100	b
13	57	b
14	100	b
15	100	b
16	100	b
17	100	b
18	98	b
19	100	b
20	100	b
21	100	b
22	100	b
23	100	b
24	100	b
25	100	b
26	51	b
27	28	b
28	100	b
29	39	b
30	57	b
31	100	b
32	100	b
33	100	b
34	100	b
35	100	b
36	100	b
37	100	b
38	100	b
39	100	b
40	100	b
41	100	b
42	100	b
43	100	b
44	100	b
45	100	b
1	35	c
2	35	c
3	35	c
4	27	c
5	37	c
6	34	c
7	27	c
8	36	c
9	31	c
10	40	c
11	43	c
12	47	c
13	34	c
14	32	c
15	44	c
16	36	c
17	50	c
18	44	c
19	39	c
20	45	c
21	55	c
22	59	c
23	52	c
24	53	c
25	100	c
26	31	c
27	37	c
28	33	c
29	40	c
30	39	c
31	34	c
32	40	c
33	33	c
34	37	c
35	46	c
36	44	c
37	100	c
38	100	c
39	43	c
40	35	c
41	40	c
42	56	c
43	100	c
44	100	c
45	100	c
1	14	d
2	15	d
3	12	d
4	13	d
5	12	d
6	100	d
7	100	d
8	13	d
9	18	d
10	100	d
11	100	d
12	49	d
13	24	d
14	100	d
15	31	d
16	33	d
17	100	d
18	32	d
19	100	d
20	28	d
21	100	d
22	100	d
23	100	d
24	100	d
25	100	d
26	13	d
27	13	d
28	11	d
29	13	d
30	27	d
31	100	d
32	100	d
33	100	d
34	100	d
35	80	d
36	100	d
37	100	d
38	100	d
39	100	d
40	100	d
41	100	d
42	100	d
43	100	d
44	100	d
45	100	d
1	17	e
2	17	e
3	14	e
4	15	e
5	14	e
6	15	e
7	15	e
8	16	e
9	17	e
10	100	e
11	60	e
12	30	e
13	31	e
14	23	e
15	42	e
16	35	e
17	100	e
18	100	e
19	100	e
20	38	e
21	71	e
22	100	e
23	58	e
24	71	e
25	100	e
26	15	e
27	17	e
28	17	e
29	18	e
30	17	e
31	25	e
32	24	e
33	18	e
34	30	e
35	40	e
36	100	e
37	52	e
38	47	e
39	100	e
40	44	e
41	100	e
42	38	e
43	100	e
44	100	e
45	100	e
};    
\addplot[line width=1pt, color = black, 
] coordinates {(1,37) (45,37)};
\addplot[line width=1pt, color = black!50, 
] coordinates {(1,84) (45,84)};
\addplot[line width=1pt, color = black!20, 
] coordinates {(1,48) (45,48)};
\addplot[line width=0.8pt, color = blue, 
] coordinates {(1,68) (45,68)};
\addplot[line width=0.8pt, color = red, 
] coordinates {(1,49) (45,49)};

\addplot[line width=1pt, color = black, 
style = dashed
] coordinates {(25.5,0) (25.5,100)};
\addplot[line width=1pt, color = black, 
style = dashed
] coordinates {(26.5,0) (26.5,100)};
\addplot[line width=1pt, color = black, 
style = dashed
] coordinates {(42.5,0) (42.5,100)};

\end{axis}
\end{tikzpicture}